\documentclass[11pt,a4paper]{article}
\usepackage{amsmath,amssymb,amsthm,mathtools}
\numberwithin{equation}{section}
\usepackage{fontspec,unicode-math}
\usepackage{longtable,booktabs,array}

\usepackage[margin=25mm]{geometry}
\usepackage{setspace}
\allowdisplaybreaks[2]
\usepackage{hyperref,bookmark}
\hypersetup{hypertexnames=false,hidelinks,pdftitle={Symphonic Joins and Finite-Product Hopf Constructions},pdflang={en-US}}
\newtheoremstyle{upright}{6pt}{6pt}{\normalfont}{0pt}{\bfseries}{.}{0.5em}{}
\theoremstyle{upright}
\newtheorem{theorem}{Theorem}
\newtheorem{lemma}{Lemma}[section]
\newtheorem{corollary}[lemma]{Corollary}
\newtheorem{proposition}[lemma]{Proposition}
\newtheorem{remark}[lemma]{Remark}

\newcommand{\R}{\mathbb R}
\newcommand{\dd}{\,dt}
\newcommand{\loc}{\mathrm{loc}}
\DeclareMathOperator{\Div}{div}
\DeclareMathOperator{\Vol}{Vol}
\DeclareMathOperator{\tr}{tr}
\DeclareMathOperator{\id}{id}

\newcommand{\vol}{\operatorname{Vol}}
\newcommand{\divg}{\operatorname{div}}
\newcommand{\Var}{\operatorname{Var}}
\newcommand{\cE}{\mathcal E}
\newcommand{\cT}{\mathcal T}
\newcommand{\cR}{\mathcal R}
\newcommand{\cJ}{\mathcal J}
\title{Symphonic Joins and Finite-Product Hopf Constructions}
\author{Xiangzhi Cao\thanks{aaa7756kijlp@163.com}\,\,\thanks{Nanjing Xiaozhuang University}}
\date{}
\begin{document}
\maketitle
\begin{abstract}
We study how collapse at the ends of a warped product affects symphonic joins and Hopf constructions. For a Hopf construction with a fixed angular map on a finite product of closed Riemannian manifolds, we characterize when the reduced energy has infimum zero among functions connecting the two target poles. Under linear collapse this occurs precisely when, at some end, the collapsing factors have total dimension greater than three, or have total dimension three and vanishing angular density on the collapsing block. In either case the infimum is unattained. For ellipsoidal joins with nonconstant spherical factor maps, smooth extension of a connecting symphonic join forces each factor map to be the restriction of a linear isometric embedding. We also derive the coupled equations for orthogonally split targets, establish a convexity criterion for minimizing separated joins, and construct smooth quaternionic examples with all angular factors active.
\end{abstract}
\newpage
\begingroup\setstretch{1}\small\tableofcontents\endgroup
\newpage

\section{Introduction}\label{sec:intro}
Symphonic maps are critical points of the integral of the squared norm of the pullback metric. The local equations allow each target component to depend on every source factor; two target blocks do not mean two source factors. The global results for the radial equation concern separated joins with factor maps $f_i:M_i\to N_i$, and Hopf constructions with one map $f:\prod_{i=1}^rP_i\to N$. The latter map may depend on all factors and has mixed angular pullback blocks. In both problems we study minimization of the reduced energy, the full map equation, and regularity across collapsed ends. Ellipsoidal geometries enter as explicit specializations rather than as restrictions on the general angular source.

Smith studied harmonic maps between spheres \cite{Smith72}. Ratto and Ding studied the harmonic Hopf construction \cite{RattoHopf,Ding94}, and Eells and Ratto treated harmonic maps between spheres and ellipsoids \cite{ER}. Further harmonic Hopf results appear in \cite{DFL,Gastel04}. For the pullback-metric functional we use the variational setting of \cite{KN,NT}. Nakauchi and Takakuwa treated the symphonic join of spherical maps \cite{NTa}. The author's earlier preprint \cite{Cao} treats joins and Hopf constructions between ellipsoids. In its Hopf construction the angular source is a product of two spheres. Here a single angular map may depend on every factor of an arbitrary finite product. Theorem~\ref{new:hopf-threshold} characterizes when connecting functions can concentrate at an end with arbitrarily small energy, and Theorem~\ref{new:linear} identifies the angular maps compatible with smooth extension of a connecting ellipsoidal join.

There is also a correction to the existence statement in \cite[Theorem 3.1]{Cao} for a prescribed angular map. Proposition~\ref{prop:slices} gives a map $S^1\times S^1\to S^3$ whose two families of slices are symphonic and whose three angular densities are constant, yet the full angular equation fails on every regular slice. Consequently this angular map admits no symphonic Hopf construction connecting the two poles. The existence results below impose the full angular equation.

Theorem~\ref{coupled:thm:full} gives the component equations. Corollaries~\ref{multi:thm:master} and~\ref{new:master} treat angular maps independent of the base. The minimizing join theorem (Theorem~\ref{new:lapses}) fixes the angular maps. Theorem~\ref{new:linear} gives a necessary condition for a connecting ellipsoidal join to extend smoothly. For Hopf constructions we distinguish existence under a trace condition (Theorem~\ref{hc:thm:exist}), the equality case of a sharp energy bound (Theorem~\ref{hc:thm:calibration}), and concentration at a collapsing end (Theorem~\ref{new:hopf-threshold}). Passing from a minimizer of the reduced energy to a symphonic map requires the angular equation and the radial equation pointwise.

Spheres, when used, carry the unit round metric, and $\Delta_g=\Div_g\nabla$. For a smooth map $f:(M,g)\to(N,g_N)$, set
\begin{equation}
S=f^*g_N,\qquad E(f)=\frac12\int_M\|S\|_g^2\,dv_g,
\qquad F(f)=2E(f).
\label{eq:energy}
\end{equation}
In a real orthonormal frame $\{e_i\}$, define
\begin{equation}\label{eq:concise:001}
\sigma_f(X)=\sum_i S(X,e_i)df(e_i).
\end{equation}
The symphonic equation is $\Div_g\sigma_f=0$, with the pullback connection.
When a map is smooth only on the regular source region, its displayed energy denotes the integral over that region. Statements about smooth symphonicity refer to that region; the global weak equation is specified below. In general, $\|f^*g_N\|^2$ is not $|df|^4$. For the variational theory and spherical joins, see \cite{KN,NT,NTa}.

\section{General geometric data and conventions}\label{g:sec:data}

Let $r,m\ge1$ be independent integers. Write $P=\prod_{j=1}^rM_j$ for a product of arbitrary Riemannian manifolds and $X=I\times P$, with $x=(x_1,\ldots,x_r)$ and full local source coordinates $q=(t,x)$. Equip $X$ with any smooth Riemannian metric $g$; a split metric will be imposed only in subsequent reductions. Let $(C,\eta)$ be a radial target manifold and $N=\prod_{i=1}^mN_i$, with a smooth family of Riemannian metrics $H_y$ on $N$, $y\in C$. On $C\times N$ put
\begin{equation}\label{coupled:eq:setup}
 G=\eta+H_y,\qquad
 F(t,x)=\bigl(u(t,x),f_1(t,x),\ldots,f_m(t,x)\bigr).
\end{equation}
Every $f_i:I\times\prod_{j=1}^rM_j\to N_i$ may depend on $t$ and on every $x_j$. We call $u(t,x)$ the radial component. When the radial target is an interval and the radial component depends only on $t$, we write $u(t,x)=\varphi(t)$. The angular maps may still depend on $t$. In particular, $m=r$ is not assumed. The source-product notation records the variables and is not needed for the local equation, which applies to any Riemannian source. On the target, the stated orthogonal splitting $G=\eta+H_y$, with $\eta$ independent of the angular coordinates, is retained.

For the equation derivations we use the auxiliary normalization $\cE=\frac{E}{2}$, where $E$ retains the convention \eqref{eq:energy}:
\begin{equation}
 \cE(F)=\frac14\int |F^*G|_g^2\mathrm d\vol_g.\label{g:eq:energy}
\end{equation}
Changing the positive prefactor does not change its critical points. If $S=F^*G$, $S^{ab}=g^{ac}g^{bd}S_{cd}$ denotes raising both indices with the source metric; it does \emph{not} denote the inverse of $S$. This notation remains valid when $F$ is not an immersion.

In source coordinates $z^a$ and target coordinates $y^A$, the symphonic tension is
\begin{equation}
 \cT^A(F)=\frac1{\sqrt{|g|}}\partial_a\bigl(\sqrt{|g|}S^{ab}\partial_bF^A\bigr)
 +S^{ab}\Gamma^A_{BC}(F)\partial_aF^B\partial_bF^C.\label{g:eq:tension}
\end{equation}
The first variation is $\delta\cE(F)[V]=-\int G(\cT(F),V)\mathrm d\vol_g$ for compactly supported variations. A smooth map is symphonic if $\cT(F)=0$. For maps of lower regularity, \emph{weakly symphonic} means that the first variation vanishes against the specified test variations.

\subsection{The full equation with independent source and target counts}
Use $\mu,\nu$ for all source coordinate indices, $A,B,D$ for radial target indices, and $I,J,K,L$ for angular target indices. Write $f=(f_1,\ldots,f_m)$ in angular coordinates. Define
\begin{equation}\label{coupled:eq:pullback}
 S_{\mu\nu}=\eta_{AB}(u)u_\mu^Au_\nu^B
                   +(H_u)_{IJ}(f)f_\mu^If_\nu^J,
 \qquad S^{\mu\nu}=g^{\mu\alpha}g^{\nu\beta}S_{\alpha\beta},
\end{equation}
where subscripts on maps denote differentiation. Let $v_g=\sqrt{\det g}$.

\begin{theorem}[Complete coupled symphonic equation]\label{coupled:thm:full}
For any smooth Riemannian source (in particular, $X=I\times P$ above) and $C^2$ component maps into the target in \eqref{coupled:eq:setup}, the full tension splits as $\mathcal T(F)=(\mathcal R,\mathcal A)$, where
\begin{align}
 \mathcal R^A={}&\frac1{v_g}\partial_\mu
           (v_gS^{\mu\nu}u_\nu^A)
       +S^{\mu\nu}\Gamma^A_{BD}(\eta)(u)u_\mu^Bu_\nu^D
       \notag\\
       &-\frac12\eta^{AB}(u)(\partial_B H_{IJ})(u,f)
                         S^{\mu\nu}f_\mu^If_\nu^J,
                \label{coupled:eq:radial}\\
 \mathcal A^I={}&\frac1{v_g}\partial_\mu
           (v_gS^{\mu\nu}f_\nu^I)
       +S^{\mu\nu}\Gamma^I_{JK}(H_u)(f)f_\mu^Jf_\nu^K
       \notag\\
       &+H^{IL}(u,f)(\partial_AH_{LJ})(u,f)
                         S^{\mu\nu}u_\mu^Af_\nu^J.
                \label{coupled:eq:angular}
\end{align}
The symbols $\Gamma(H_y)$ use angular differentiation at fixed $y$; $\partial_AH$ differentiates the radial parameter at fixed angular coordinates. The matrix $H^{IL}$ is the inverse angular target metric, whereas $S^{\mu\nu}$ denotes raising both indices, as above. Symphonicity is equivalent to $\mathcal R^A=0$ and $\mathcal A^I=0$ for every component, pointwise. All sums include $t$ and every angular source coordinate. No compactness, curvature or separation hypothesis on the maps is required.
\end{theorem}
\begin{proof}
For the orthogonal target metric, its nonzero connection blocks, apart from the symmetries in the lower indices, are
\[
 \Gamma^A_{BD}=\Gamma^A_{BD}(\eta),\qquad
 \Gamma^I_{JK}=\Gamma^I_{JK}(H_y),\qquad
 \Gamma^A_{IJ}=-\tfrac12\eta^{AB}\partial_BH_{IJ},\qquad
 \Gamma^I_{AJ}=\tfrac12H^{IL}\partial_AH_{LJ}.
\]
The other mixed blocks $\Gamma^A_{BI}$ and $\Gamma^I_{AB}$ vanish. Substitute these identities and \eqref{coupled:eq:pullback} into \eqref{g:eq:tension}. The two equal mixed contributions to the angular equation combine, giving the coefficient one in the last term of \eqref{coupled:eq:angular}. The two components together give the first-variation equation for arbitrary compactly supported target variations.
\end{proof}

\begin{corollary}[Explicit equations for $m$ evolving angular maps]\label{coupled:cor:warped}
Suppose $H_y=\sum_{i=1}^m b_i(y)^2\gamma_i$, with $b_i>0$. Set $K_i=f_i^*\gamma_i$, as a tensor on the full source $X$, and $K_0=u^*\eta$. Then $S=K_0+\sum_i b_i(u)^2K_i$. For a map $w$ into a fixed metric target, write
\[
 \mathscr L_S w=\frac1{v_g}\nabla^w_\mu
                 (v_gS^{\mu\nu}\partial_\nu w),
 \qquad \langle S,K_i\rangle_g=S^{\mu\nu}(K_i)_{\mu\nu},
\]
where $\nabla^w$ acts on the target index and uses $\eta$ for $u$ and $\gamma_i$ for $f_i$. The complete system is
\begin{align}
 0&=\mathscr L_Su
       -\sum_{i=1}^m b_i(u)\operatorname{grad}_{\eta}b_i(u)
                     \langle S,K_i\rangle_g,
               \label{coupled:eq:warped-radial}\\
 0&=\mathscr L_Sf_i
       +2S^{\mu\nu}\partial_\mu(\log b_i\circ u)\partial_\nu f_i,
       \qquad 1\le i\le m.
               \label{coupled:eq:warped-angular}
\end{align}
Equivalently, the second equation is
\begin{equation}\label{coupled:eq:weighted-angular}
 \frac1{v_g b_i(u)^2}\nabla^{f_i}_\mu
       \bigl(v_g b_i(u)^2S^{\mu\nu}\partial_\nu f_i\bigr)=0.
\end{equation}
\end{corollary}
\begin{proof}
In Theorem~\ref{coupled:thm:full}, the radial derivative of the $i$th metric block is $2b_i\partial_Ab_i\,\gamma_i$ and $H_i^{-1}\partial_AH_i=2\partial_A\log b_i$ times the identity. Substitution gives both equations. The product rule proves \eqref{coupled:eq:weighted-angular}.
\end{proof}

For the split source metric $g=h(t)^2dt^2+g_t$, keep $u=u(t,x)$ and $f_i=f_i(t,x)$ fully coupled. With $a,b$ denoting angular source coordinate indices, the pullback blocks are
\begin{align}
 S_{00}&=\eta(u_t,u_t)+\sum_{i=1}^m b_i(u)^2\gamma_i(\partial_tf_i,\partial_tf_i),
         \label{coupled:eq:blocks}\\
 S_{0a}&=\eta(u_t,u_a)+\sum_{i=1}^m b_i(u)^2\gamma_i(\partial_tf_i,\partial_af_i),
         \notag\\
 S_{ab}&=\eta(u_a,u_b)+\sum_{i=1}^m b_i(u)^2\gamma_i(\partial_af_i,\partial_bf_i).
         \notag
\end{align}
Here $u_t=\partial_tu$ and $u_a=\partial_au$ are tangent vectors to $C$ along $u$. For a scalar radial coordinate with $\eta=k(s)^2ds^2$, the three radial contributions are respectively $k(u)^2u_t^2$, $k(u)^2u_tu_a$ and $k(u)^2u_au_b$.
With $\theta=S_{0a}dx^a$ and $W=S_{ab}dx^a dx^b$, the energy density is still
\begin{equation}\label{coupled:eq:energy-blocks}
 |S|_g^2=h^{-4}S_{00}^2+2h^{-2}|\theta|_{g_t}^2+|W|_{g_t}^2.
\end{equation}

For this one-dimensional radial target, the radial equation \eqref{coupled:eq:warped-radial} takes the explicit full-source form
\begin{equation}\label{general-u:eq:scalar-pde}
 \begin{split}
 0={}&\frac1{v_g}\partial_\mu
          \bigl(v_g S^{\mu\nu}\partial_\nu u\bigr)
       +\frac{k'(u)}{k(u)}S^{\mu\nu}\partial_\mu u\partial_\nu u\\
     &-\sum_{i=1}^m\frac{b_i(u)b_i'(u)}{k(u)^2}
                         \langle S,K_i\rangle_g.
 \end{split}
\end{equation}
It is accompanied by all $m$ angular equations \eqref{coupled:eq:warped-angular}, with
$\partial_\mu(\log b_i\circ u)=\frac{b_i'(u)}{b_i(u)}\partial_\mu u$.
Every $\mu,\nu$ sum includes $t$ and all angular source coordinates, so mixed and angular second derivatives are retained. Formula \eqref{general-u:eq:scalar-pde} follows directly from Corollary~\ref{coupled:cor:warped}, since $\Gamma^s_{ss}=\frac{k'}{k}$ and $\operatorname{grad}_{\eta}b_i=\frac{b_i'}{k^2}\partial_s$.

The specialization order is explicit. First, $u(t,x)=\varphi(t)$ removes $u_a$ but leaves all $\partial_tf_i$ terms. Next, imposing $\partial_tf_i=0$ for every $i$ gives the separated reductions below. Neither condition alone justifies discarding every mixed pullback block. In particular, if the $f_i$ are independent of $t$ but $u_a\ne0$, the contribution $\eta(u_t,u_a)$ remains. The counts $r$ and $m$ stay independent throughout; selecting a two-block target never forces $r=2$. The complete system is used only to derive these specializations, and no general existence result for $u(t,x)$ is asserted.

\subsection{Independent factor data for later specializations}
For the following independent-factor data only, impose $m=r$ and let each $f_i$ depend solely on $x_i$. This is a specialization of \eqref{coupled:eq:setup}, not a restriction on the complete system.

Let $I,J$ be open intervals, and let $(M_i,g_i)$ and $(N_i,\gamma_i)$ be connected Riemannian manifolds, $1\leq i\leq r$, with $m_i=\dim M_i$. On their products consider
\begin{align}
 g&=h(t)^2\mathrm d t^2+\sum_{i=1}^r a_i(t)^2g_i,
 &G&=k(s)^2\mathrm d s^2+\sum_{i=1}^r b_i(s)^2\gamma_i,\label{g:eq:metrics}
\end{align}
where all the displayed scalar functions are smooth and strictly positive on their intervals. Completeness, constant curvature, homogeneity, and embeddings into Euclidean space are not assumed. Write
\[
 \rho(t)=h(t)\prod_i a_i(t)^{m_i}.
\]
For smooth maps $f_i:M_i\to N_i$ and any $\varphi\in C^2(I,J)$, set
\begin{equation}
 F(t,x_1,\ldots,x_r)=(\varphi(t),f_1(x_1),\ldots,f_r(x_r)).\label{g:eq:ansatz}
\end{equation}
For each factor define
\begin{equation}
 S_i=f_i^*\gamma_i,\qquad q_i=|S_i|_{g_i}^2,\qquad
 T_i=\divg_{g_i}\sigma_{f_i}.\label{g:eq:factor}
\end{equation}
Here $\sigma_{f_i}(X)=\sum_\alpha\gamma_i(df_i(X),df_i(e_\alpha))df_i(e_\alpha)$ in a $g_i$-orthonormal frame. Thus $q_i$ is the unnormalized quartic density, not the ordinary Dirichlet density $|df_i|^2$.

\subsection{Base-independent angular reduction}
We now specialize Theorem~\ref{coupled:thm:full} by fixing the angular maps along the source base. Source and target base dimensions may still be arbitrary.

Let $(B,\beta)$ and $(C,\eta)$ be Riemannian manifolds, with local coordinates $p^\alpha$ and $y^A$. Let $g_p$ be a smooth family of metrics on an arbitrary manifold $P$, and let $H_y$ be a smooth family of metrics on $N=\prod_{i=1}^mN_i$, where $m\ge1$ is finite. Each $f_i:P\to N_i$ may depend on all angular variables. Put
\begin{equation}\label{multi:eq:setup}
 g=\beta+g_p,\qquad G=\eta+H_y,\qquad
 F(p,z)=\bigl(u(p),f_1(z),\ldots,f_m(z)\bigr),
\end{equation}
where $u:B\to C$ and $f=(f_1,\ldots,f_m):P\to N$ are $C^2$. The metrics $H_y$ need not split into target factor metrics. The absence of base-angular metric blocks and the independence of the $f_i$ from $p$ are retained assumptions.

\begin{corollary}[Base-independent angular maps]\label{multi:thm:master}
Set $A=u^*\eta$, raise its two indices with $\beta$, and define
\[
 \mu(p,z)=\sqrt{\det\beta(p)}\sqrt{\det g_p(z)},\qquad
 Q(p,y,z)=|f^*H_y|_{g_p}^2.
\]
The complete symphonic tension of \eqref{multi:eq:setup} is the pair
\begin{align}
 \mathcal T(F)&=\bigl(\mathcal R,\mathcal T_{g_p,H_{u(p)}}(f)\bigr),
       \label{multi:eq:master}\\
 \mathcal R^A
 &=\frac1\mu\partial_\alpha
       \bigl(\mu A^{\alpha\beta}\partial_\beta u^A\bigr)
   +A^{\alpha\beta}\Gamma^A_{DE}(\eta)(u)
       \partial_\alpha u^D\partial_\beta u^E
   -\frac14(\operatorname{grad}_{\eta}Q)^A(p,u,z).
       \label{multi:eq:base}
\end{align}
The gradient differentiates only $y$, holding $p,z$ fixed. Thus symphonicity is equivalent to the vanishing of both components at every $(p,z)$, without angular averaging. This is a local equivalence, with no compactness or curvature hypothesis.
\end{corollary}
\begin{proof}
Specialize Theorem~\ref{coupled:thm:full} to $u=u(p)$, $f=f(z)$ and $g=\beta+g_p$. Writing $W=f^*H_y$, the pullback blocks become $S_{\alpha\beta}=A_{\alpha\beta}$, $S_{ab}=W_{ab}$ and $S_{\alpha a}=0$. The radial angular-metric term is $-\tfrac14\operatorname{grad}_{\eta}Q$, since $\partial_DQ=2W^{ab}(\partial_DH_{IJ})(f)\partial_af^I\partial_bf^J$. All mixed terms in \eqref{coupled:eq:angular} vanish. Its remaining divergence is the slice tension, because the base volume factor is angularly constant. This gives both displayed components.
\end{proof}

\begin{corollary}[Explicit coupled equations for multiple target factors]\label{multi:cor:warped}
Suppose $H_y=\sum_{i=1}^m b_i(y)^2\gamma_i$, with $b_i>0$. For the maps $f_i:P\to(N_i,\gamma_i)$ define
\[
 H_i=f_i^*\gamma_i,\qquad C_{ij}(p,z)=\langle H_i,H_j\rangle_{g_p}.
\]
For a symmetric covariant tensor $K$ on $P$, define the section along $f_i$
\begin{equation}\label{multi:eq:angular-operator}
 \mathscr D_i[K]
 =\frac1{\sqrt{\det g_p}}\nabla^{f_i}_a
       \left(\sqrt{\det g_p}\,K^{ab}\partial_bf_i\right),
\end{equation}
where both indices of $K$ are raised with $g_p$ and $\nabla^{f_i}$ uses $\gamma_i$. Then
\begin{align}
 Q&=\sum_{i,j=1}^m b_i^2b_j^2 C_{ij},
       \label{multi:eq:cross-energy}\\
 \frac14\operatorname{grad}_{\eta}Q
   &=\sum_{i=1}^m b_i\left(\sum_{j=1}^m b_j^2C_{ij}\right)
                   \operatorname{grad}_{\eta}b_i,
       \label{multi:eq:cross-force}\\
 \bigl(\mathcal T_{g_p,H_u}(f)\bigr)_i
   &=\sum_{j=1}^m b_j(u)^2\mathscr D_i[H_j].
       \label{multi:eq:cross-angular}
\end{align}
Together with \eqref{multi:eq:base}, these are necessary and sufficient pointwise equations for $F$ to be symphonic. In particular, when the $f_i$ share angular variables, symphonicity of each individual $f_i$ does not by itself remove the $j\ne i$ terms.
\end{corollary}
\begin{proof}
The angular pullback is $W=\sum_jb_j^2H_j$. Expanding its squared norm gives \eqref{multi:eq:cross-energy}; differentiating and using $C_{ij}=C_{ji}$ gives \eqref{multi:eq:cross-force}. At fixed $p$, each $b_j(u)$ is angularly constant. The angular target connection is the product connection of the $\gamma_i$, so its $i$th component is the operator \eqref{multi:eq:angular-operator} applied to $W$. This proves \eqref{multi:eq:cross-angular}.
\end{proof}
If $m=r$, $P=\prod_iM_i$, $g_p=\sum_i a_i(p)^2g_i$ and $f_i=f_i(x_i)$, then $C_{ij}=0$ for $i\ne j$, $C_{ii}=a_i^{-4}q_i$, and $\mathscr D_i[H_j]=0$ for $i\ne j$. These conclusions use the factor dependence as well as the orthogonal source metric; they cannot be inferred merely from having several target components.

\subsection{One-dimensional specialization with arbitrary slice metrics}
The product factors in \eqref{g:eq:metrics} are useful for construction, but are not needed for the local reduction itself. The following statement also permits mixed angular blocks and angular dependence of the volume density.

\begin{corollary}[General slice-metric reduction]\label{new:master}
Let $P,N$ be arbitrary smooth manifolds, let $g_t$ and $H_s$ be smooth families of Riemannian metrics on them, and let $h(t),k(s)>0$. On $I\times P$ and $J\times N$ put
\[
 g=h(t)^2dt^2+g_t,\qquad G=k(s)^2ds^2+H_s.
\]
For $F(t,z)=(\varphi(t),f(z))$, with $f:P\to N$ smooth and $\varphi\in C^2(I,J)$, define in local angular coordinates
\[
 \mu(t,z)=h(t)\sqrt{\det g_t(z)},\qquad
 Q(t,s,z)=|f^*H_s|_{g_t}^{2}.
\]
Then the full symphonic tension, in the normalization \eqref{g:eq:tension}, is
\begin{align}
 \mathcal T(F)&=k(\varphi)^{-2}\mathcal R(t,z)\partial_s
                 +\mathcal T_{g_t,H_{\varphi(t)}}(f),\label{new:eq:master}\\
 \mathcal R(t,z)&=\frac1\mu\partial_t\left(\mu\frac{k(\varphi)^4}{h^4}\varphi'^3\right)
 -\frac{k(\varphi)^3k'(\varphi)}{h^4}\varphi'^4
 -\frac14\partial_sQ(t,\varphi(t),z).\label{new:eq:master-radial}
\end{align}
In the last term $\partial_s$ holds $t,z$ fixed. Thus $F$ is symphonic precisely when both displayed components vanish pointwise. The statement is local, does not require compactness, and imposes no curvature condition on either angular manifold.
\end{corollary}
\begin{proof}
In Corollary~\ref{multi:thm:master}, take $B=I$, $\beta=h^2dt^2$, $C=J$, $\eta=k^2ds^2$, $u=\varphi$, and one angular target $N$. Then $A^{00}=\frac{k^2\varphi'^2}{h^4}$ and $\Gamma^s_{ss}=\frac{k'}{k}$. Thus
\[
 \mathcal T^s=\mu^{-1}\partial_t(\mu k^2h^{-4}\varphi'^3)
                  +kk'h^{-4}\varphi'^4-(4k^2)^{-1}Q_s.
\]
Multiplying by $k^2$ and applying the product rule gives \eqref{new:eq:master-radial}; the angular component is already supplied by that corollary.
\end{proof}
The coordinate density $\mu$ changes by a $t$-independent angular Jacobian under a change of angular chart, so the expression is intrinsic. In contrast to a multiply warped product, $\partial_t\log\mu$ can depend on $z$. That dependence is an additional compatibility condition, and cannot be replaced by an angular average. This corollary retains an orthogonal interval splitting and a $t$-independent angular map. Radial-angular cross terms are outside this separated ansatz.

\section{Full reduction and compatibility}
\begin{corollary}[Full separated reduction]\label{g:thm:reduction}
For the metrics \eqref{g:eq:metrics} and map \eqref{g:eq:ansatz},
\begin{equation}
 \cT(F)=k(\varphi)^{-2}\cR(t,x)\,\partial_s
       +\sum_i\frac{b_i(\varphi)^2}{a_i(t)^4}T_i,\label{g:eq:split}
\end{equation}
where each $T_i$ is inserted into its target factor and
\begin{align}
 \cR(t,x)={}&\frac1\rho\frac{\mathrm d}{\mathrm d t}
 \left(\rho\frac{k(\varphi)^4}{h^4}(\varphi')^3\right)
 -\frac{k(\varphi)^3k'(\varphi)}{h^4}(\varphi')^4\notag\\
 &-\sum_i\frac{b_i(\varphi)^3b_i'(\varphi)}{a_i^4}q_i(x_i).
 \label{g:eq:radial}
\end{align}
Primes on $k,b_i$ mean differentiation in $s$, while primes on $\varphi$ mean differentiation in $t$. Consequently, $F$ is symphonic if and only if every $f_i$ is symphonic and $\cR(t,x)=0$ for every $(t,x)$.
\end{corollary}
\begin{proof}
Apply Corollary~\ref{multi:cor:warped} with $B=I$, $\beta=h^2dt^2$, $C=J$, $\eta=k^2ds^2$, $g_t=\sum a_i^2g_i$ and $f_i=f_i(x_i)$. Its cross terms vanish, giving $Q=\sum b_i^4a_i^{-4}q_i$ and angular components $b_i^2a_i^{-4}T_i$. The one-dimensional base conversion in Corollary~\ref{new:master} gives the displayed radial formula. Positivity of the scale factors gives the equivalence.
\end{proof}

\begin{corollary}[Necessary constant densities]\label{g:cor:constant}
Suppose $F$ is symphonic. If for some $t_i\in I$,
\begin{equation}
 b_i'(\varphi(t_i))\neq0,\label{g:eq:detected}
\end{equation}
then $q_i$ is constant on $M_i$. If \eqref{g:eq:detected} holds for every $i$, symphonicity is equivalent to the following three conditions: $T_i=0$, $q_i\equiv\lambda_i$ for every factor, and \eqref{g:eq:radial} with $q_i$ replaced by $\lambda_i$.
\end{corollary}
\begin{proof}
Fix $t=t_i$ and all angular variables except $x_i$. Subtract the radial equation at two points of $M_i$. The coefficient $\frac{b_i^3b_i'}{a_i^4}$ is nonzero, so the two values of $q_i$ agree. Sufficiency follows from Corollary~\ref{g:thm:reduction}.
\end{proof}
\begin{remark}
No independence condition on different warping functions is needed: the factor variables are independent. The exceptional case $b_i'(\varphi(t))\equiv0$ is real. The radial equation then contains no $q_i$-term, so it imposes no constant-density requirement on that factor. Angular symphonicity remains necessary.
\end{remark}

\subsection{An exact residual identity}
Assume here that all $M_i$ are closed, and integrate over their product using the product of normalized volume measures, denoted by $\mathrm d\mu$. Define
\[
 \bar q_i=\int_{M_i}q_i\,\frac{\mathrm d\vol_{g_i}}{\vol(M_i)},\qquad
 \Var(q_i)=\int_{M_i}(q_i-\bar q_i)^2\,\frac{\mathrm d\vol_{g_i}}{\vol(M_i)}.
\]
Let $c_i(t)=\frac{b_i(\varphi)^3b_i'(\varphi)}{a_i^4}$, and let $\overline{\cR}(t)$ denote \eqref{g:eq:radial} with $q_i$ replaced by $\bar q_i$.
\begin{theorem}[Variance obstruction]\label{g:thm:variance}
For every smooth $\varphi:I\to J$,
\begin{equation}
 \int\cR(t,x)^2\mathrm d\mu(x)=\overline{\cR}(t)^2+\sum_i c_i(t)^2\Var(q_i).\label{g:eq:variance}
\end{equation}
In particular,
\begin{equation}
 \int|\cT(F)|_G^2\mathrm d\mu\ \geq\
 k(\varphi)^{-2}\left(\overline{\cR}^2+\sum_i c_i^2\Var(q_i)\right),\label{g:eq:residualbound}
\end{equation}
with equality when every $f_i$ is symphonic. If $\overline{\cR}=0$, the full radial equation vanishes exactly when $c_i(t)^2\Var(q_i)=0$ for every $i,t$.
\end{theorem}
\begin{proof}
Subtract the averaged equation:
\[
 \cR=\overline{\cR}-\sum_i c_i(q_i-\bar q_i).
\]
Every centered term has zero mean. For $i\neq j$, independence of the product variables gives
$\int(q_i-\bar q_i)(q_j-\bar q_j)\mathrm d\mu=0$. Expanding the square proves \eqref{g:eq:variance}. The radial vector is $k^{-2}\cR\partial_s$, whose squared target norm is $k^{-2}\cR^2$. Orthogonality to the angular factors proves \eqref{g:eq:residualbound}.
\end{proof}
Even when the averaged radial residual vanishes, a nonzero variance term gives a positive lower bound for the mean squared tension.

\subsection{Independent factor maps with several radial variables}
In this specialization $m=r$ and $f_i=f_i(x_i)$; the general coupled equations above allow independent $r,m$.
\begin{corollary}[Weighted base equation]\label{multi:cor:independent}
Let $a_i:B\to(0,\infty)$, $b_i:C\to(0,\infty)$, and
\[
 g=\beta+\sum_i a_i^2g_i,\qquad
 G=\eta+\sum_i b_i^2\gamma_i,\qquad
 F=(u(p),f_1(x_1),\ldots,f_r(x_r)).
\]
Put $\omega=\prod_i a_i^{m_i}$ and let
$\sigma_u(X)=\sum_\alpha\eta(du(X),du(e_\alpha))du(e_\alpha)$ for a $\beta$-orthonormal frame. With $q_i,T_i$ as in \eqref{g:eq:factor}, the full equations are
\begin{equation}\label{multi:eq:weighted-base}
 T_i=0\quad(1\le i\le r),\qquad
 \frac1\omega\operatorname{div}_{\beta}^{u}(\omega\sigma_u)
 =\sum_{i=1}^r\frac{q_i(x_i)}{a_i^4}
                     b_i(u)^3\operatorname{grad}_{\eta}b_i(u).
\end{equation}
The divergence uses the pullback connection along $u$. If
$\operatorname{grad}_{\eta}b_i(u(p_i))\ne0$ at some $p_i$, then a symphonic $F$ forces $q_i$ to be constant on $M_i$.
\end{corollary}
\begin{proof}
Apply Corollary~\ref{multi:cor:warped}. The cross terms vanish as stated there, the angular component is $b_i^2a_i^{-4}T_i$, and the base volume density factors as $\sqrt{\det\beta}\,\omega$ times a $p$-independent angular density. These facts give \eqref{multi:eq:weighted-base}. For the last assertion, fix $p_i$ and every angular variable except $x_i$, and subtract the vector equation at two choices of $x_i$.
\end{proof}
This statement allows both base dimensions to exceed one. In particular, it applies on the regular part of full multi-factor joins when the base and warping functions have the prescribed join geometry. It does not assert smoothness of their collapsed boundaries.

\subsection{The complete multi-factor ellipsoid equation}
Let the target block count be $m\ge2$ and $c_i>0$. Define the radial manifold and its induced Euclidean metric by
\[
 C_c=\left\{y\in(0,\infty)^m:\sum_i \frac{y_i^2}{c_i^2}=1\right\},
 \qquad \eta=\left.\sum_i dy_i^2\right|_{C_c}.
\]
For unit round spheres $S^{n_i}$, define
\[
 \iota(y,\xi_1,\ldots,\xi_m)=(y_1\xi_1,\ldots,y_m\xi_m).
\]
This map is an isometry from $(C_c\times\prod_iS^{n_i},\eta+\sum_i y_i^2g_{S^{n_i}})$ onto the open subset with all blocks nonzero of the ellipsoid
\begin{equation}\label{multi:eq:ellipsoid-target}
 E_c=\left\{(Y_1,\ldots,Y_m):
        \sum_i\frac{|Y_i|^2}{c_i^2}=1\right\}.
\end{equation}
The radial manifold of this ellipsoid has dimension $m-1$.

Indeed, $|\xi_i|=1$ and $\langle\xi_i,d\xi_i\rangle=0$ give
\begin{equation}\label{coupled:eq:ellipsoid-isometry}
 \iota^*g_{E_c}
   =\left.\sum_i dy_i^2\right|_{TC_c}
                 +\sum_i y_i^2g_{S^{n_i}}
   =\eta+\sum_i y_i^2g_{S^{n_i}}.
\end{equation}
The inverse is $y_i=|Y_i|$, $\xi_i=\frac{Y_i}{|Y_i|}$ on this open set. Thus $\iota$ identifies the entire regular target, rather than merely embedding a lower-dimensional slice. If
\[
 F(t,x)=(y(t,x),f_1(t,x),\ldots,f_m(t,x)),\qquad \Phi=\iota\circ F,
\]
then preservation of the metric and Levi--Civita connection gives
\begin{equation}\label{coupled:eq:tension-isometry}
 \mathcal T_{E_c}(\Phi)=d\iota\bigl(\mathcal T_{\eta+\sum y_i^2g_{S^{n_i}}}(F)\bigr).
\end{equation}
Here $x=(x_1,\ldots,x_r)$ may contain any number of source factors, independently of $m$. Since $d\iota$ is invertible, the two full map equations are equivalent. An arbitrary isometric immersion into a larger target would not justify this conclusion; here the isometry is onto the regular part of $E_c$ with its induced metric.

\begin{corollary}[Paired factor maps into an ellipsoid]\label{multi:cor:ellipsoid}
In this paired specialization only, set $m=r$ and let $f_i$ depend on $x_i$ alone. The two-target-block case with arbitrary source count is instead covered by Corollary~\ref{blocks:cor:fixed}.
Let $g=h(t)^2dt^2+\sum_i a_i(t)^2g_i$ on $I\times\prod_iM_i$, where the $M_i$ are arbitrary connected Riemannian manifolds. Let $f_i:M_i\to S^{n_i}$ be smooth, $y:I\to C_c$ be $C^2$, and set
\begin{equation}\label{multi:eq:ellipsoid-map}
 \Phi(t,x_1,\ldots,x_r)
       =\bigl(y_1(t)f_1(x_1),\ldots,y_r(t)f_r(x_r)\bigr).
\end{equation}
Write $\rho=h\prod_i a_i^{m_i}$, $v^2=|y'|^2$, and $q_i=|f_i^*g_{S^{n_i}}|_{g_i}^2$. Then $\Phi$ is symphonic if and only if every $f_i$ is symphonic, each $q_i$ is a constant $\lambda_i\ge0$, and
\begin{equation}\label{multi:eq:ellipsoid-system}
 \frac1\rho\frac d{dt}\left(\frac{\rho}{h^4}|y'|^2y'\right)
 -\left(\frac{\lambda_i y_i^3}{a_i^4}\right)_{i=1}^r
 =\Lambda(t)\left(\frac{y_i}{c_i^2}\right)_{i=1}^r,
 \qquad \sum_i\frac{y_i^2}{c_i^2}=1.
\end{equation}
The multiplier is determined, rather than an additional free datum:
\begin{equation}\label{multi:eq:multiplier}
 \Lambda
 =-\frac{\displaystyle
      \frac{|y'|^2}{h^4}\sum_i\frac{(y_i')^2}{c_i^2}
       +\sum_i\frac{\lambda_i y_i^4}{a_i^4c_i^2}}
       {\displaystyle\sum_i\frac{y_i^2}{c_i^4}}.
\end{equation}
All assertions concern the regular region $y_i>0$.
\end{corollary}
\begin{proof}
First use \eqref{coupled:eq:ellipsoid-isometry}--\eqref{coupled:eq:tension-isometry} to identify $\Phi$ with the intrinsic map $F$. Apply Corollary~\ref{multi:cor:independent} to this $F$, taking $B=I$, $C=C_c$, $u=y$ and $b_i(y)=y_i$. On $C_c$, each $\operatorname{grad}_{\eta}y_i$ is nonzero: the vector $e_i$ cannot be parallel to the normal $n=(\frac{y_j}{c_j^2})_j$, which has every entry positive and $r\ge2$. Thus all $q_i$ are constant. The intrinsic base equation is the vanishing of the tangent projection of the left side of \eqref{multi:eq:ellipsoid-system}. Its orthogonal complement is the line spanned by $n$, giving that system. Differentiating the constraint yields $\langle y',n\rangle=0$ and $\langle y'',n\rangle=-\sum_i\frac{(y_i')^2}{c_i^2}$. Taking the normal inner product gives \eqref{multi:eq:multiplier}. Conversely, the system has zero tangent projection and the factor tensions vanish, so the full equation follows from the same corollary.
\end{proof}
For closed $M_i$ with the stated constant densities, division of $\mathcal E$ in \eqref{g:eq:energy} by $\prod_i\operatorname{Vol}(M_i)$ reduces the energy to
\begin{equation}\label{multi:eq:ellipsoid-energy}
 J[y]=\frac14\int_I\rho(t)
       \left(\frac{|y'|^4}{h^4}
                +\sum_i\frac{\lambda_i y_i^4}{a_i^4}\right)dt,
       \qquad y(t)\in C_c.
\end{equation}
Its constrained Euler equation is \eqref{multi:eq:ellipsoid-system}. For the full map equation the factor maps must also be symphonic.

\begin{corollary}[Exact condition for restriction to a radial curve]\label{multi:cor:curve}
Under the constant-density factor hypotheses above, let $R:J\to C_c$ be a smooth unit-speed curve, $y(t)=R(s(t))$, and set
\[
 V_t(y)=\frac14\sum_i\frac{\lambda_i y_i^4}{a_i(t)^4},\qquad
 \kappa=\nabla^{C_c}_{R'}R'.
\]
Then the full radial system is equivalent to the pair
\begin{align}
 \frac1\rho\left(\frac\rho{h^4}s'^3\right)'
       &=\frac d{ds}V_t(R(s)),\label{multi:eq:curve-tangent}\\
 \frac{s'^4}{h^4}\kappa(s)
       &=\bigl(\operatorname{grad}_{\eta}V_t(R(s))\bigr)^{\perp R'}.
       \label{multi:eq:curve-normal}
\end{align}
In the first line $t$ is held fixed in the derivative on the right. The second line is the component normal to the curve within $C_c$.
\end{corollary}
\begin{proof}
Since $|R'|=1$, the covariant derivative of $\rho h^{-4}s'^3R'$ divided by $\rho$ is
$\rho^{-1}(\rho h^{-4}s'^3)'R'+h^{-4}s'^4\kappa$.
Project the intrinsic equation onto $R'$ and its orthogonal complement in $TC_c$.
\end{proof}
The radial dimension is determined by the target count $m$. When $m=2$, there is no transverse radial direction. In the paired subclass above, where $r=m$, taking $y_1=c\cos\varphi$, $y_2=d\sin\varphi$ recovers \eqref{g:eq:radial}, with $b_1=c\cos s$, $b_2=d\sin s$ and $k(s)^2=c^2\sin^2s+d^2\cos^2s$. When $m\ge3$, \eqref{multi:eq:curve-normal} has $m-2$ transverse components; solving only \eqref{multi:eq:curve-tangent} is not sufficient. These reductions explain exactly when a map whose radial component lies on a prescribed curve is symphonic as a map into the full ellipsoid.

\subsection{Two target blocks with arbitrary source-factor count}\label{blocks:sec:two-target}
Here ``two target blocks'' refers to the two displayed Euclidean components of the target map. It imposes no restriction on the number $r$ of source angular factors. Write $P=\prod_{j=1}^rM_j$ and $x=(x_1,\ldots,x_r)$, with arbitrary finite $r\ge1$. The symbol $m$ in the general equation counts angular maps, which must be distinguished from the number of displayed ambient blocks: the Hopf form has one angular map and one radial scalar block.
\begin{center}
\begin{tabular}{@{}lll@{}}
\toprule
Target form & Source count & Angular maps in the full system\\
\midrule
Two angular target blocks & arbitrary $r$ & $m=2$: $f_1(t,x),f_2(t,x)$\\
Hopf-type two-block form & arbitrary $r$ & $m=1$: $f(t,x)$\\
\bottomrule
\end{tabular}
\end{center}

For two angular target blocks, let $f_i:I\times P\to S^{n_i}$, $i=1,2$, and let $u:I\times P\to(0,\frac{\pi}{2})$ be $C^2$. The map is
\begin{equation}\label{blocks:eq:evolving-join}
 \Phi(t,x)=\bigl(c\cos u(t,x)f_1(t,x),\ d\sin u(t,x)f_2(t,x)\bigr).
\end{equation}
Take $C=(0,\frac{\pi}{2})$, $\eta=k(s)^2ds^2$, $b_1=c\cos s$, $b_2=d\sin s$ and $k^2=c^2\sin^2s+d^2\cos^2s$. The isometry
\[
 \iota_J(s,\xi_1,\xi_2)=(c\cos s\,\xi_1,d\sin s\,\xi_2)
\]
identifies \eqref{blocks:eq:evolving-join} with $F=(u(t,x),f_1(t,x),f_2(t,x))$. Its complete equation is precisely Corollary~\ref{coupled:cor:warped} with $m=2$ and arbitrary $r$. In particular, the angular derivatives of $u$, the $t$-derivatives of both $f_i$, and the full radial-angular blocks in \eqref{coupled:eq:blocks} remain present.

\begin{corollary}[Two angular target blocks on an arbitrary finite product]\label{blocks:cor:fixed}
Let $r\ge1$ be arbitrary, let
$g=h(t)^2dt^2+g_t$, $g_t=\sum_{j=1}^r a_j(t)^2g_j$ on $I\times P$, with smooth positive $h,a_j$, and set $\rho=h\prod_{j=1}^r a_j^{d_j}$, $d_j=\dim M_j$. Impose the additional specialization $u(t,x)=\varphi(t)$, with $\varphi\in C^2(I,(0,\frac{\pi}{2}))$. Suppose the two smooth angular maps are independent of $t$: $f_i:P\to S^{n_i}$, $i=1,2$. They may each depend on all $r$ source factors. Set
\[
 H_i=f_i^*g_{S^{n_i}},\qquad C_{ij}(t,x)=\langle H_i,H_j\rangle_{g_t},
\]
and
\begin{equation}\label{blocks:eq:cross-density}
 Q(t,s,x)=c^4\cos^4s\,C_{11}+d^4\sin^4s\,C_{22}
                  +2c^2d^2\cos^2s\sin^2s\,C_{12}.
\end{equation}
Then \eqref{blocks:eq:evolving-join}, with $f_i(t,x)=f_i(x)$, is symphonic on its regular region if and only if
\begin{align}
 0={}&\frac1\rho\left(\rho\frac{k(\varphi)^4}{h^4}\varphi'^3\right)'
        -\frac{k(\varphi)^3k'(\varphi)}{h^4}\varphi'^4
        -\frac14\partial_s Q(t,\varphi(t),x),
        \label{blocks:eq:radial}\\
 0={}&c^2\cos^2\varphi\,\mathscr D_i[H_1]
            +d^2\sin^2\varphi\,\mathscr D_i[H_2],
               \qquad i=1,2,
        \label{blocks:eq:angular}
\end{align}
pointwise for every $(t,x)$, where $\mathscr D_i$ is the slice operator in \eqref{multi:eq:angular-operator}. In $\partial_sQ$, both $t$ and $x$ are held fixed.
\end{corollary}
\begin{proof}
Use Corollary~\ref{multi:cor:warped} with two angular targets and $P=\prod_{j=1}^rM_j$. Expanding $|b_1^2H_1+b_2^2H_2|_{g_t}^2$ gives \eqref{blocks:eq:cross-density}. Corollary~\ref{new:master} gives its radial equation, and \eqref{multi:eq:cross-angular} gives the two angular equations. The isometry $\iota_J$ transfers the full equation to the ellipsoid, as in \eqref{coupled:eq:tension-isometry}.
\end{proof}
Neither $C_{12}$ nor the cross operators $\mathscr D_i[H_j]$ for $i\ne j$ can be omitted for general maps on the common source $P$. Individual symphonicity of the $f_i$, or separate constancy of $C_{11}$ and $C_{22}$, is not asserted by this corollary. Under the further specialization $r=2$, $f_1=f_1(x_1)$ and $f_2=f_2(x_2)$, these cross terms vanish and the paired equation \eqref{g:eq:radial} follows. They can also vanish for larger $r$, for example when the two maps depend on disjoint groups of orthogonal source factors. The later two-source-sphere results remain a narrower special case, not the definition of a two-block target.

For the Hopf-type two-block form, with the same arbitrary $P$, take
\begin{equation}\label{blocks:eq:evolving-hopf}
 \Psi(t,x)=\bigl(c\cos u(t,x)f(t,x),\ d\sin u(t,x)\bigr),
 \qquad f:I\times P\to S^{n-1}.
\end{equation}
Here $C=(-\frac{\pi}{2},\frac{\pi}{2})$, $\eta=k(s)^2ds^2$, $b_1=c\cos s$ and $m=1$. The isometry
\[
 \iota_H(s,\xi)=(c\cos s\,\xi,d\sin s)
\]
gives its full equation directly from Corollary~\ref{coupled:cor:warped}. The second displayed target block is the radial scalar $d\sin u(t,x)$, not a second freely varying angular map. When $u(t,x)=\varphi(t)$ and $\partial_tf=0$, this specializes to the Hopf construction with $f=f(x_1,\ldots,x_r)$: Corollary~\ref{new:master} yields \eqref{hc:eq:angular}--\eqref{hc:eq:radial} after target arclength change and \eqref{eq:hopf-angular}--\eqref{eq:hopf-radial} in ellipsoid latitude. All mixed entries among the source factors are retained. These local deductions impose no restriction $r=2$; the global completion and existence statements retain their individual hypotheses.

\section{Mixed angular maps: exact rank and covariance criteria}\label{new:sec:rank}
Let $P=\prod_{i=1}^r M_i$, $g_t=\sum_i a_i(t)^2g_i$, and let $f:P\to(N,\gamma)$ depend on all factors. The target metric is $k(s)^2ds^2+b(s)^2\gamma$. In fixed $g_i$-orthonormal frames set
\[
 q_{ii}=\sum_{\alpha,\beta}\langle df(e_{i,\alpha}),df(e_{i,\beta})\rangle^2,
 \qquad q_{ij}=\sum_{\alpha,\beta}\langle df(e_{i,\alpha}),df(e_{j,\beta})\rangle^2\quad(i<j).
\]
Let $\mathbf q$ list these $d=\frac{r(r+1)}{2}$ quantities, and let $\mathbf w(t)$ list $a_i^{-4}$ in the diagonal positions and $2a_i^{-2}a_j^{-2}$ in the other positions. Then
\begin{equation}\label{new:eq:mixed-density}
 Q_t(z)=|f^*\gamma|_{g_t}^{2}=\mathbf w(t)^T\mathbf q(z).
\end{equation}
Corollary~\ref{new:master} gives the angular equation $b(\varphi)^2\mathcal T_{g_t,\gamma}(f)=0$ and the radial equation
\begin{equation}\label{new:eq:mixed-radial}
 D(t)-b(\varphi)^3b'(\varphi)\mathbf w(t)^T\mathbf q(z)=0,
\quad D(t)=\rho^{-1}(\rho k^4h^{-4}\varphi'^3)'-k^3k'h^{-4}\varphi'^4,
\end{equation}
where $\rho=h\prod_i a_i^{m_i}$. No hypothesis about the separate slices of $f$ replaces the displayed full angular equation.

\begin{theorem}[Precisely the detected combinations are constant]\label{new:rank}
Suppose the full map is symphonic, and let $J_0\subset I$ be an open interval on which $b(\varphi)^3b'(\varphi)\ne0$. Put $W=\operatorname{span}\{\mathbf w(t):t\in J_0\}$. Then
\begin{equation}\label{new:eq:kernel}
 \mathbf q(z)-\mathbf q(z_0)\in W^\perp\qquad(z,z_0\in P).
\end{equation}
In particular, full span $W=\mathbb R^d$ forces every diagonal and mixed density to be constant. Conversely, if \eqref{new:eq:kernel} holds and the radial equation holds at one angular point for all $t\in J_0$, it holds at all angular points there. This converse concerns the radial equation; the angular equation is still required.

If the factors are closed, let $d\nu_P$ be normalized product volume and let $\Sigma=\int_P(\mathbf q-\overline{\mathbf q})(\mathbf q-\overline{\mathbf q})^T d\nu_P$. For any finite positive measure $d\nu$ on a compact subinterval of $J_0$, put $\mathsf G=\int\mathbf w\mathbf w^T d\nu$. Then
\begin{equation}\label{new:eq:gram}
 \int\operatorname{Var}_{\nu_P}(Q_t)\,d\nu(t)=\tr(\mathsf G\Sigma).
\end{equation}
If $\mathsf G$ is positive definite, it follows that
\begin{equation}\label{new:eq:quant-rank}
 \sum_{i\le j}\operatorname{Var}_{\nu_P}(q_{ij})
 \le\lambda_{\min}(\mathsf G)^{-1}\int\operatorname{Var}_{\nu_P}(Q_t)\,d\nu(t).
\end{equation}
\end{theorem}
\begin{proof}
Subtract \eqref{new:eq:mixed-radial} at $z,z_0$ and divide by the nonzero common coefficient. The resulting orthogonality for every $t$ is exactly \eqref{new:eq:kernel}; the converse follows by adding it back. For the quantitative statement, \eqref{new:eq:mixed-density} gives $\operatorname{Var}(Q_t)=\mathbf w(t)^T\Sigma\mathbf w(t)$. Integration gives \eqref{new:eq:gram}. Since $\Sigma$ is positive semidefinite, $\tr(\mathsf G\Sigma)\ge\lambda_{\min}(\mathsf G)\tr\Sigma$, which proves the estimate.
\end{proof}
Writing $\overline{\mathcal R}$ for the average of the radial residual gives the further exact identity
\[
 \int_P\mathcal R(t,z)^2d\nu_P(z)=\overline{\mathcal R}(t)^2+
 b(\varphi)^6b'(\varphi)^2\operatorname{Var}_{\nu_P}(Q_t).
\]
Multiplication by $k(\varphi)^{-2}$ is a lower bound for the averaged squared norm of the full tension, with equality when the angular tension vanishes. Unlike the separated variance identity, mixed densities need not be statistically independent; the covariance matrix is essential. Rank deficiency specifies which combinations the radial equation fails to detect. It does not assert the existence of a map realizing arbitrary values in that kernel.

For the two ellipsoidal source warps the three weights are linearly independent on every open subinterval. Indeed, after division by $a_2(t)^{-4}$ they are $x(t)^2,1,2x(t)$, where $x=\frac{a_2(t)^2}{a_1(t)^2}=\frac{b^2}{a^2}\tan^2t$ is nonconstant. A polynomial of degree at most two vanishing on an interval is zero. This recovers the three-density rigidity of Corollary~\ref{prop:hopf-rigidity}.

\section{An exact hidden-convexity criterion}\label{new:sec:convex}
To express the reduced energy in a general radial coordinate, let $v=\chi(s)$ be any strictly increasing smooth coordinate with $\chi'>0$ in the open target interval. For the separated energy define
\[
 \beta(v)=\left(\frac{\chi'(s)}{k(s)}\right)^{\frac{4}{3}},\qquad A(t)=\frac{\rho(t)}{h(t)^4},
 \qquad s=\chi^{-1}(v),
\]
\[
 U(t,v)=\frac{\rho(t)}4\sum_i\frac{\lambda_i b_i(\chi^{-1}(v))^4}{a_i(t)^4}.
\]
The integrand of the reduced energy is
\begin{equation}\label{new:eq:convex-density}
 \mathscr L(t,v,z)=\frac{A(t)z^4}{4\beta(v)^3}+U(t,v),\qquad z=v'.
\end{equation}
\begin{theorem}[Necessary and sufficient convexity in a chosen coordinate]\label{new:convex}
Assume $A(t)>0$, $\beta>0$ and $\beta,U(t,\cdot)$ are $C^2$ on an open interval. For fixed $t$, $\mathscr L(t,\cdot,\cdot)$ is jointly convex for all $z\in\mathbb R$ if and only if
\begin{equation}\label{new:eq:convex-iff}
 \beta''\le0\quad\hbox{and}\quad U_{vv}(t,\cdot)\ge0.
\end{equation}
If $\beta''\le0$ and $U_{vv}(t,v)\ge\lambda(t)\ge0$, and a minimizing $v_*$ and a competitor $v$ belong to an admissible class closed under affine interpolation in $v$, then
\begin{equation}\label{new:eq:general-gap}
 \mathcal J(v)-\mathcal J(v_*)\ge\frac12\int\lambda(t)|v-v_*|^2\,dt,
\end{equation}
whenever the energies are finite. In particular $\lambda>0$ almost everywhere gives uniqueness in this class.
\end{theorem}
\begin{proof}
Differentiation gives
\[
 \mathscr L_{zz}=\frac{3Az^2}{\beta^3},\quad
 \mathscr L_{vv}=U_{vv}+\frac{3Az^4}{4\beta^5}(4\beta'^2-\beta\beta''),
\]
\begin{equation}\label{new:eq:hessian-det}
 \det D^2_{(v,z)}\mathscr L
 =\frac{3Az^2}{\beta^3}U_{vv}-\frac{9A^2z^6}{4\beta^7}\beta''.
\end{equation}
Under \eqref{new:eq:convex-iff}, the Hessian is positive semidefinite, including at $z=0$, where it is $\operatorname{diag}(U_{vv},0)$. Necessity of $U_{vv}\ge0$ follows at $z=0$; a positive $\beta''$ would make the determinant negative for large $|z|$. This proves the equivalence. Convexity of the derivative part and strong convexity of the potential give
\[
 \mathcal J((1-\theta)v_*+\theta v)\le(1-\theta)\mathcal J(v_*)+\theta\mathcal J(v)
 -\frac{\theta(1-\theta)}2\int\lambda|v-v_*|^2.
\]
Use minimality, divide by $\theta$, and let $\theta\downarrow0$ to obtain \eqref{new:eq:general-gap}. Truncation of the nonnegative gap integral justifies the same argument if its finiteness was not known in advance.
\end{proof}
This is a sharp criterion for this scalar affine-interpolation mechanism. It is not a necessary condition for existence or uniqueness by other methods, and it does not ensure existence of a suitable coordinate $\chi$. The perspective argument and hidden-convexity principle have established precedents in \cite{BF}. For ellipsoids $\chi(s)=\sin^2s$ yields a positive constant multiple of the function $\beta$ in \eqref{eq:convex-variables}, so the criterion includes the all-dimension join theorem rather than imposing convexity in target arclength.

\section{Target arclength and the variational problem}
Introduce target arclength
\[
 u=K(s)=\int_{s_*}^{s}k(\xi)\mathrm d\xi,\qquad B_i(u)=b_i(K^{-1}(u)).
\]
The source coordinate is left unchanged. Set
\begin{equation}
 A(t)=\frac{\rho(t)}{h(t)^4},\qquad
 P(t,u)=\frac{\rho(t)}4\sum_i\frac{\lambda_i B_i(u)^4}{a_i(t)^4}.\label{g:eq:AP}
\end{equation}
If $q_i\equiv\lambda_i$, integration over closed factors gives
\begin{equation}
 \frac{\cE(F)}{\prod_i\vol(M_i)}=\cJ[u]
 :=\int_I\left(\frac14 A(t)|u'|^4+P(t,u)\right)\mathrm d t.\label{g:eq:J}
\end{equation}
Indeed, the block decomposition yields
\[
 |F^*G|_g^2=\frac{k(\varphi)^4}{h^4}(\varphi')^4
     +\sum_i\frac{b_i(\varphi)^4}{a_i^4}q_i.
\]
Since $u'=k(\varphi)\varphi'$, no additional radial factor remains in \eqref{g:eq:J}. The radial equation becomes
\begin{equation}
 (A(u')^3)'=P_u(t,u).\label{g:eq:flux}
\end{equation}
For nonconstant densities the same integrated energy contains $\bar q_i$ instead; Theorem~\ref{g:thm:variance} explains why its critical points need not be symphonic maps.

\subsection{Degenerate weights with controlled traces}
\begin{theorem}[Weighted existence and convergence]\label{g:thm:weighted}
Let $A>0$ be continuous on $(0,L)$ and assume
\begin{equation}
 H:=\int_0^L A(t)^{-\frac{1}{3}}\mathrm d t<\infty.\label{g:eq:tracecondition}
\end{equation}
Let $P:(0,L)\times\mathbb R\to[0,\infty)$ be measurable in $t$ and continuous in $u$ for almost every $t$. Consider locally absolutely continuous functions with finite $\int A|u'|^4$ and fixed endpoint traces. If at least one such function has finite $\cJ$, then $\cJ$ has a minimizer $u_*$. If, in addition, $P(t,\cdot)$ is convex almost everywhere, this minimizer is unique and, for every admissible finite-energy $v$,
\begin{equation}
 \cJ[v]-\cJ[u_*]\geq\frac1{32}\int_0^L A|v'-u_*'|^4\mathrm d t.\label{g:eq:weightedgap}
\end{equation}
Under the convexity hypothesis every minimizing sequence converges uniformly on $[0,L]$ and strongly in its weighted derivative norm. If $P$ is $C^1$ on compact interior rectangles and $A$ is locally $C^1$, the minimizer solves \eqref{g:eq:flux} in the interior and is locally $C^{1,\frac{1}{3}}$.
\end{theorem}
\begin{proof}
Weighted H\"older gives, for $0<s<t<L$,
\begin{equation}
 |u(t)-u(s)|\leq
 \left(\int_s^t A|u'|^4\mathrm d\tau\right)^{\frac{1}{4}}
 \left(\int_s^t A^{-\frac{1}{3}}\mathrm d\tau\right)^{\frac{3}{4}}.\label{g:eq:holder}
\end{equation}
Condition \eqref{g:eq:tracecondition} ensures endpoint traces and a common modulus of continuity for an energy-bounded sequence with one fixed trace. Thus a minimizing sequence has a uniformly convergent subsequence on $[0,L]$. Its derivatives converge weakly in $L^4(A\mathrm d t)$; testing against compactly supported smooth functions identifies the limit as the derivative of the uniform limit. Weak lower semicontinuity handles the kinetic term, and Fatou's lemma handles the potential. The traces persist, proving existence. If $P$ is convex, strict convexity of the derivative term and the fixed traces give uniqueness.

For real $a,b$,
\[
 \frac12\left(\frac{a^4}4+\frac{b^4}4\right)
 -\frac14\left(\frac{a+b}2\right)^4\geq\frac{|a-b|^4}{64}.
\]
Apply this and convexity of $P$ to $\frac{u_*+v}{2}$. Minimality implies
\[
 \cJ[u_*]\leq\cJ[\frac{u_*+v}{2}]
 \leq\tfrac12\cJ[u_*]+\tfrac12\cJ[v]-\tfrac1{64}\int A|u_*'-v'|^4,
\]
which proves \eqref{g:eq:weightedgap} without requiring an integrable global first variation of the possibly singular potential. The gap and \eqref{g:eq:holder}, applied to $v-u_*$, give convergence. Compactly supported variations and the flux argument prove the interior assertion.
\end{proof}
If an additional weight $\rho$ is integrable, uniform convergence also gives strong convergence in $L^4(\rho\mathrm d t)$. If $h$ and $h^{-1}$ are bounded, then $A=\frac{\rho}{h^4}$ is comparable to $\rho$, so the derivative convergence can also be stated in the $\rho$-weighted norm. These comparisons must be checked before changing weighted spaces.

\subsection{Regular intervals}
\begin{corollary}[Existence, uniqueness, and an energy gap]\label{g:thm:regular}
Let $I=(0,L)$, let $A\in C^1([0,L])$ with $0<A_*\leq A\leq A^*$, and let $P\in C^1([0,L]\times\R)$ be nonnegative. For prescribed real endpoint values $u_-,u_+$, the functional \eqref{g:eq:J} has a minimizer in
\[
 \mathcal A=\{u\in W^{1,4}(0,L):u(0)=u_-,\ u(L)=u_+\}.
\]
If $P(t,\cdot)$ is convex for every $t$, the minimizer $u_*$ is unique and
\begin{equation}
 \cJ[v]-\cJ[u_*]\geq\frac1{12}\int_0^L A|v'-u_*'|^4\mathrm d t
 \quad(v\in\mathcal A).\label{g:eq:gap}
\end{equation}
Moreover $u_*\in C^{1,\frac{1}{3}}([0,L])$, its flux $A(u_*')^3$ is $C^1$, and it solves \eqref{g:eq:flux} weakly. For smooth coefficients it is smooth on the open set where $u_*'\neq0$. Under the convexity hypothesis every minimizing sequence converges strongly in $W^{1,4}$ and uniformly on $[0,L]$.
\end{corollary}
\begin{proof}
The weight satisfies $\int A^{-\frac{1}{3}}<\infty$, and the affine competitor has finite energy. Theorem~\ref{g:thm:weighted} therefore gives existence, and under convexity uniqueness and convergence. Since $A$ is bounded above and below, weighted derivative convergence and the fixed traces give strong $W^{1,4}$ convergence. The flux equation has continuous right side up to the endpoints, so its flux is $C^1$ there and the cube-root formula yields $C^{1,\frac{1}{3}}$ regularity; away from zero derivative the smooth equation is nondegenerate.
For the improved constant, use
\begin{equation}\label{g:eq:quartic}
 \tfrac14 b^4-\tfrac14 a^4-a^3(b-a)\ge\tfrac1{12}(b-a)^4.
\end{equation}
The difference between the two sides is $\tfrac32(b-a)^2(a+\frac{b-a}{3})^2$. Integrate this inequality and the supporting-line inequality for $P$, then subtract the first variation tested with $v-u_*$. All terms are integrable on this regular interval, giving \eqref{g:eq:gap}.
\end{proof}
\begin{remark}
Under the convexity hypothesis, the minimizer of the reduced energy is unique with the angular maps fixed. It does not prove uniqueness among all maps with the same boundary data, nor stability against arbitrary angular variations. Also, a $C^{1,\frac{1}{3}}$ solution of the radial equation is not automatically smooth at a zero of $u'$. At this regularity Corollary~\ref{g:thm:reduction} applies in divergence form and yields a weakly symphonic map when its angular conditions hold.
\end{remark}

A simple sufficient geometric condition for convexity is
\begin{equation}
 \frac{\mathrm d^2}{\mathrm d u^2}B_i(u)^4\geq0
 \quad\text{for each factor with }\lambda_i>0.\label{g:eq:convex}
\end{equation}
This is a substantive hypothesis on the target. The theorem does not assert it for spherical or ellipsoidal warping functions on their whole joining intervals. For existence without convexity the minimizer need not be unique.

When applying this abstract arclength theorem to a target chart, one must have $K(J)=\mathbb R$ or separately verify that the minimizer stays in $K(J)$. Extending the scalar potential beyond its geometric interval does not by itself establish that condition. The all-dimension theorem below instead treats the prescribed compact interval for the values of $u$.

\section{Examples with general factors}
\subsection{A nonlinear Dirichlet class}
Take arbitrary closed factors $(M_i,g_i)$, put $(N_i,\gamma_i)=(M_i,g_i)$, and let $f_i$ be the identity. Then $\nabla df_i=0$, $T_i=0$, and $q_i=m_i$. Choose any smooth positive source functions $h,a_i$ on $[0,L]$ and, in target arclength, choose
\[
 B_i(u)=\exp(\kappa_i u),\qquad u\in\R.
\]
Since $(B_i^4)''=16\kappa_i^2e^{4\kappa_i u}\geq0$, Corollary~\ref{g:thm:regular} gives a unique minimizer for any two prescribed endpoint values. The associated separated map is weakly symphonic for variations with fixed boundary values, and is smooth where its radial derivative is nonzero. The factors may have variable curvature and need not be spheres. 

\subsection{An explicit smooth family on nonisometric products}\label{g:subsec:seed}
Let $h=1$, let all $a_i$ be positive constants, and let all target $B_i$ be positive constants. For arbitrary symphonic factor maps and endpoint values $u_-\neq u_+$,
\begin{equation}
 F_0(t,x)=\left(u_-+\frac{u_+-u_-}{L}t,f_1(x_1),\ldots,f_r(x_r)\right)\label{g:eq:explicit}
\end{equation}
is smooth and symphonic. Indeed the flux is constant and all radial angular coefficients vanish. No constant-density assumption on the factor maps is needed in this exceptional product case. If the $f_i$ are immersions, then $F_0$ is an immersion. Choosing different constant scales gives source and target metrics that need not be isometric.

Taking the identity maps of arbitrary closed factors gives explicit separated cylinder examples with curved factors as well as flat tori.

\section{Endpoint capacity and smooth collapse}\label{g:sec:endpoints}
\subsection{The quartic endpoint threshold}
\begin{proposition}[Exact transition cost]\label{g:prop:capacity}
For $0<\delta<\epsilon<L$ and a positive continuous weight $A$, the infimum over absolutely continuous functions on $[\delta,\epsilon]$ with the indicated traces is
\begin{equation}
 \inf_{\substack{v(\delta)=1\\v(\epsilon)=0}}
 \int_\delta^\epsilon A|v'|^4\mathrm d t
 =\left(\int_\delta^\epsilon A^{-\frac{1}{3}}\mathrm d t\right)^{-3}.\label{g:eq:capacity}
\end{equation}
If $A(t)\asymp t^\beta$ near $0$, the derivative energy controls an endpoint trace through \eqref{g:eq:holder} when $\beta<3$. For $\beta\geq3$, the endpoint has zero capacity for this derivative energy: cutoffs equal to $1$ near $0$ and $0$ beyond any fixed $\epsilon$ can have arbitrarily small derivative energy.
\end{proposition}
\begin{proof}
H\"older applied to $1=|\int_\delta^\epsilon v'|$ gives the lower bound. Equality is attained by
\[
 v(t)=\frac{\int_t^\epsilon A(\tau)^{-\frac{1}{3}}\mathrm d\tau}
                {\int_\delta^\epsilon A(\tau)^{-\frac{1}{3}}\mathrm d\tau}.
\]
Extend this function by $1$ to the left of $\delta$ and by $0$ to the right of $\epsilon$. The trace integral converges at zero exactly when $\beta<3$. If it diverges, sending $\delta$ to zero makes the displayed energy tend to zero, including logarithmically at $\beta=3$.
\end{proof}
\begin{remark}
This is a statement about the derivative part of the energy. A singular angular potential can impose additional endpoint conditions, so zero derivative capacity alone is not a nonexistence result for the full join equation. Conversely, finite weighted energy without further hypotheses does not justify retention of prescribed boundary values in a minimizing limit.
\end{remark}
For source arclength $r$, a collapsing factor of dimension $m$ with scale $a(r)\sim r$, all other scales bounded away from zero, gives $A(r)\asymp r^m$. The threshold is therefore $m=3$. This helps explain why regular-interval variational arguments cannot be transferred unchanged to every collapsed join.

\subsection{Arbitrary links need not give a smooth manifold}
\begin{proposition}[Smooth-link obstruction]\label{g:prop:link}
Suppose a metric near a putative smooth collapsed stratum has the form
\[
 \mathrm d r^2+a(r)^2g_L+g_B(r),\qquad \frac{a(r)}{r}\longrightarrow1,
\]
where $r$ is normal distance and the connected compact factor $L$ collapses to the center of each normal slice. If the metric extends smoothly and nondegenerately across that stratum as a Riemannian manifold, then $(L,g_L)$ is isometric to the unit round sphere of the normal dimension minus one.
\end{proposition}
\begin{proof}
Fix a point of the stratum and restrict to a normal slice. Under rescaling by $r^{-2}$, the metric link of the slice converges to $(L,g_L)$ by the assumed asymptotic form. At a smooth Riemannian point the tangent cone is Euclidean, and its metric unit link is the unit round sphere. Uniqueness of that tangent metric forces the stated isometry. In dimension two this also fixes the circumference of the circle link, ruling out a nontrivial cone angle.
\end{proof}
For a leading coefficient $\frac{a(r)}{r}\to c>0$, the same statement concerns $(L,c^2g_L)$. A spherical space-form quotient can instead give an orbifold cone; it is not a smooth-manifold endpoint of this type. Additional parity and expansion conditions on the warping functions are needed even when the link is round.

The map has its own extension condition. For instance, when both source and target factors collapse, one must analyze the regularity of the corresponding angular map in Cartesian normal coordinates. The prescribed limit of the radial component gives its endpoint value. Smooth extension also depends on its rate of approach and on the angular map. Thus the interior theory, weighted trace theory, and smooth geometric extension are three separate steps.

\section{All-dimension minimizing joins beyond ellipsoids}\label{new:sec:lapses}
The following theorem proves the global scalar and weak-map conclusions before any ellipsoidal specialization. No endpoint traces are imposed on its minimizing class.

\begin{theorem}[Independent radial geometries and arbitrary angular factors]\label{new:lapses}
Let $(M_i,g_i)$ and $(N_i,\gamma_i)$ be closed connected Riemannian manifolds of positive dimensions, and let $f_i:M_i\to N_i$ be smooth symphonic maps with constant densities $|f_i^*\gamma_i|^2=\lambda_i>0$. Write $m_i=\dim M_i$. Fix $a,b,c,d>0$, $L=\frac{\pi}{2}$, $I=(0,L)$, and put
\begin{align*}
 g&=h(t)^2dt^2+a^2\cos^2t\,g_1+b^2\sin^2t\,g_2,\\
 G&=k(s)^2ds^2+c^2\cos^2s\,\gamma_1+d^2\sin^2s\,\gamma_2.
\end{align*}
Assume $h,k$ are smooth and strictly positive on $[0,L]$, and that $k$ has a smooth even extension across each endpoint. Set
\begin{equation}\label{new:eq:beta-k}
 \beta_k(v)=\left(\frac{v(1-v)}{k(\arcsin\sqrt v)^2}\right)^{\frac{2}{3}},\qquad 0<v<1.
\end{equation}
If $\beta_k''\le0$, define $\rho=h\cos^{m_1}t\sin^{m_2}t$, $a_1=\frac{c^4\lambda_1}{a^4}$, $a_2=\frac{d^4\lambda_2}{b^4}$, and define $J,X,X_0$ as follows.
We minimize the extended real-valued functional
\begin{equation}
J(\varphi)=\frac14\int_I\left\{
\frac{k(\varphi)^4}{h^4}(\varphi')^4
+a_1\frac{\cos^4\varphi}{\cos^4t}
+a_2\frac{\sin^4\varphi}{\sin^4t}\right\}\rho\dd.
\label{eq:J}
\end{equation}
The admissible space carries no prescribed endpoint traces:
\begin{equation}
X=\left\{\alpha\in W^{1,4}_{\loc}(I):\|\alpha\|_X^4
:=\int_I\rho(|\alpha|^4+|\alpha'|^4)\dd<\infty\right\},
\quad X_0=\{\alpha\in X:0\le\alpha\le L\}.
\label{eq:X}
\end{equation}
All interval constraints are understood almost everywhere.

Then:
\begin{enumerate}
\item $J$ has a unique finite-energy minimizer in the class $X_0$ without imposed endpoint traces. Its representative is continuous on $[0,L]$, smooth and strictly increasing on $(0,L)$, and connects $0$ to $L$.
\item Its arclength flux $P=\left(\frac{\rho}{h^4}\right)(K(\varphi)')^3$, where $K(s)=\int_0^s k(r)dr$, tends to zero at both endpoints. The energy-gap estimate \eqref{eq:main-energy-gap} holds, and every minimizing sequence converges strongly in the weighted space $X$.
\item The separated map is symphonic on the regular region, has finite energy, and extends continuously to the metric join completions. If these completions are smooth compact Riemannian manifolds with the usual collapsed strata, the extension is globally $W^{1,4}$ and weakly symphonic for arbitrary target variations.
\end{enumerate}
No relation between $h$ and $a,b$, or between $k$ and $c,d$, is required for the scalar assertions. Smooth completion is an additional geometric hypothesis for the final assertion.
\end{theorem}
\begin{proof}
Put $K(s)=\int_0^s k(r)\,dr$, $\Theta=K^{-1}$, $U=K(L)$ and $u=K(\varphi)$.
Even reflection extends $k$ to a positive smooth even $\pi$-periodic function on $\mathbb R$; thus $K$ is a smooth bi-Lipschitz diffeomorphism. Write
\begin{equation}\label{eq:transformed}
 \widetilde J(u)=J(\Theta(u))=\frac14\int_I(q|u'|^4+\rho W(t,u))dt,
 \qquad q=\frac{\rho}{h^4},
\end{equation}
where $A=\frac{a_1}{\cos^4t}$, $B=\frac{a_2}{\sin^4t}$ and
$W(t,u)=A\cos^4\Theta(u)+B\sin^4\Theta(u)$.
The endpoint coefficient comparisons use positivity and smoothness of $h,k$ on the closed interval.

\emph{Existence and the unconstrained flux equation.}
The graph $\alpha\mapsto(\alpha,\alpha')$ identifies $X$ with a closed subspace of $L^4(\rho\,dt)^2$: convergence there identifies the derivative by testing on compact interior intervals. Hence $X$ is reflexive and its interval-constrained subset is closed and convex. The competitor $\varphi=t$ has finite energy. Since $q\asymp\rho$ and $0\le u\le U$, a minimizing sequence is bounded in $X$. Weak convergence in $X$, local uniform convergence, weak lower semicontinuity of the derivative term and Fatou's lemma give a minimizer.

The potential $W$ is invariant under reflection about $0,U$ and translation by $2U$. The folding map $R(v)=\min_{j\in\mathbb Z}|v-2jU|$ preserves $W$ and $|v'|$ almost everywhere. Thus $R(u+\varepsilon\eta)$ is an admissible variation for every compactly supported smooth $\eta$. Differentiation gives
\begin{equation}\label{eq:flux}
 P=q(u')^3,\qquad
 P'=\frac\rho4W_u(t,u)
 =\frac{\rho\sin\varphi\cos\varphi}{k(\varphi)}
       (B\sin^2\varphi-A\cos^2\varphi).
\end{equation}
Consequently $P\in C^1(I)$ and $u,\varphi\in C^{1,\frac{1}{3}}_{\mathrm{loc}}(I)$.

\emph{Moving well and exclusion of the obstacles.}
For $0<u<U$, direct differentiation gives
\begin{equation}\label{short:eq:well}
 \begin{gathered}
 u_*(t)=K\!\left(\arctan\!\left(\sqrt{\frac{a_1}{a_2}}\tan^2t\right)\right),\qquad
 \operatorname{sgn}W_u=\operatorname{sgn}(u-u_*),\\
 W_{uu}(t,u_*)=\frac{8AB}{(A+B)k(\Theta(u_*))^2}>0.
 \end{gathered}
\end{equation}
The well $u_*$ is strictly increasing from $0$ to $U$.
If $u(t_0)=0$, then $P(t_0)=0$ and $P'\le0$ on a neighborhood where $u<u_*$. Thus $u'\ge0$ to the left and $u'\le0$ to the right; nonnegativity forces $u=0$ on that neighborhood. But
$W(t,s)=A-\frac{2As^2}{k(0)^2}+O(s^4)$, so a nonzero nonnegative bump $\eta$ supported there gives
\[
 \widetilde J(u+\varepsilon\eta)-\widetilde J(u)
 =-\frac{\varepsilon^2}{2}\int\frac{\rho A}{k(0)^2}\eta^2dt+O(\varepsilon^4)<0.
\]
This is impossible. The expansion $W(t,U-s)=B-\frac{2Bs^2}{k(L)^2}+O(s^4)$ proves the upper-obstacle assertion. Hence $0<u<U$ in $I$.

\emph{Automatic endpoint values in every dimension.}
On $(0,t_0)$ truncate $u$ down to $\max\{u(t_0),u_*(t_0)\}$; on $(t_0,L)$ truncate it up to $\min\{u(t_0),u_*(t_0)\}$. The derivative energy cannot increase, and a nontrivial truncation strictly lowers the potential by \eqref{short:eq:well}. Therefore
\begin{equation}\label{short:eq:truncation}
 u(t)\le\max\{u(t_0),u_*(t_0)\}\ (t<t_0),\qquad
 u(t)\ge\min\{u(t_0),u_*(t_0)\}\ (t>t_0).
\end{equation}
Taking $t_0$ along sequences realizing the left liminf and right limsup proves that both endpoint limits exist. Fixing an interior $t_0$ also shows $0\le\ell_0<U$ and $0<\ell_L\le U$.

Suppose $0<\ell_0<U$ and put $m=m_2$, $h_0=h(0)$. Then
$\rho\sim h_0t^m$, $q\sim h_0^{-3}t^m$ and $B\sim a_2t^{-4}$.
For $m\le3$ the potential is nonintegrable. For $m>3$, \eqref{eq:flux} gives $P'\sim\gamma t^{m-4}$ with $\gamma>0$, so $P$ has a finite limit. A nonzero limit would make
$q|u'|^4=|P|^{\frac{4}{3}}q^{-\frac{1}{3}}\asymp t^{-\frac{m}{3}}$ nonintegrable. Thus $P(0+)=0$, whence $P\sim\frac{\gamma t^{m-3}}{m-3}$ and $u'\sim \frac{C}{t}$, contradicting boundedness. Therefore $\ell_0=0$. Reflecting $t\mapsto L-t$, $\varphi\mapsto L-\varphi$, and exchanging the two factors proves $\ell_L=U$.

\emph{Strict monotonicity and zero endpoint flux.}
If $u$ decreases, its endpoint values and strict interior bounds provide a local maximum at $r$ followed by a local minimum at $v>r$, with $u(r)>u(v)$. Their flux signs imply $P'(r)\le0$, $P'(v)\ge0$, hence
$u(r)\le u_*(r)<u_*(v)\le u(v)$, a contradiction. Thus $u'\ge0$ and $P\ge0$.
At a zero $t_0$ of $P$, also $P'(t_0)=0$, so $u=u_*$ and $u'=0$ there. Differentiating \eqref{eq:flux} and the identity $W_u(t,u_*(t))=0$ gives
\[
 P''(t_0)=-\tfrac14\rho(t_0)W_{uu}(t_0,u_*(t_0))u_*'(t_0)<0,
\]
contrary to $P\ge0$. Thus $P>0$, and the smooth system $u'=\left(\frac{P}{q}\right)^{\frac{1}{3}}$, $P'=\frac{\rho W_u}{4}$ gives interior smoothness.

Near zero, \eqref{eq:flux} gives $P'\ge-C\rho u$. The negative part of $P'$ is integrable; $P\ge0$ then has a finite nonnegative left limit. If it were positive, $u'\sim Ct^{-\frac{m}{3}}$. For $m\ge3$ this contradicts the finite trace. For $m=1,2$, integration gives $u\sim C't^{1-\frac{m}{3}}$; substituting back into \eqref{eq:flux} yields $P'\sim \frac{C''}{t}$ with $C''>0$, again contradicting the finite limit. Reflection proves
\begin{equation}\label{eq:zero-end-flux}
 P(0+)=P(L-)=0.
\end{equation}

\emph{Convexity, the gap, and strong convergence.}
Put $v=\sin^2\varphi$. The integrand of $J$ becomes
\begin{equation}\label{eq:convex-variables}
 \frac{\rho|v'|^4}{64h^4\beta_k(v)^3}
 +\frac\rho4\{A(1-v)^2+Bv^2\}.
\end{equation}
The perspective $\frac{|z|^4}{\eta^3}$ is jointly convex for $\eta>0$ and nonincreasing in $\eta$. Concavity of $\beta_k$ therefore makes the first term convex in $(v,v')$, as also follows from Theorem~\ref{new:convex}. At $v=0,1$ use the lower semicontinuous extension, with value zero for $v'=0$ and infinity otherwise. Absolutely continuous functions have derivative zero almost everywhere on their level sets. Let $v_\theta=(1-\theta)\sin^2\widehat\varphi+\theta\sin^2\psi$, $0<\theta<1$. The strict interior bounds for $\widehat\varphi$ keep $v_\theta$ uniformly away from both $0$ and $1$ on each compact subinterval of $I$. Thus $\arcsin\sqrt{v_\theta}$ is locally Sobolev by the chain rule. Convexity bounds its kinetic integral on all of $I$, and the bounds for $h,k$ give membership in $X_0$. Apply minimality at $v_\theta$, divide the resulting inequality by $\theta$, and let $\theta\downarrow0$. The exact quadratic identity for the potential and minimality give
\begin{equation}\label{eq:main-energy-gap}
 J(\psi)-J(\widehat\varphi)\ge\frac14\int_I\rho(A+B)
 (\sin^2\psi-\sin^2\widehat\varphi)^2dt.
\end{equation}
This proves uniqueness. The elementary inequality
$|\sin^2\alpha-\sin^2\gamma|\ge\frac{4|\alpha-\gamma|^2}{\pi^2}$ on $[0,L]$ also gives
\[
 \int_I\rho|\psi-\widehat\varphi|^4dt
 \le\frac{\pi^4}{4(a_1+a_2)}[J(\psi)-J(\widehat\varphi)].
\]
For a minimizing sequence set $u_j=K(\psi_j)$. Any subsequence has a further subsequence converging weakly in $X$ and locally uniformly to $\widehat u$. The kinetic and potential integrals are separately lower semicontinuous; convergence of their sum forces convergence of the kinetic integrals. Uniform convexity of $L^4(q\,dt)$ gives $u_j'\to\widehat u'$ strongly there. Finally
\[
 \psi_j'-\widehat\varphi'
 =\frac{u_j'-\widehat u'}{k(\psi_j)}
 +\widehat u'\left(\frac1{k(\psi_j)}-\frac1{k(\widehat\varphi)}\right).
\]
Boundedness of $k^{-1}$, dominated convergence and $q\asymp\rho$ give strong derivative convergence in $L^4(\rho\,dt)$. The preceding gap controls the function part. Uniqueness and the subsequence criterion give strong $X$ convergence and local uniform convergence for the whole minimizing sequence.

\emph{The complete map and its weak extension.}
Corollary~\ref{g:thm:reduction} gives the full equation on the regular region. The endpoint limits are independent of the variables collapsed there, so the map descends continuously to the metric completions. Suppose these completions are smooth. Embed the compact target isometrically into Euclidean space. Finite energy and
$|dF|^4\le(\dim M)|F^*G|^2$ give an $L^4$ differential. The boundary areas of truncated cylinders are $O(t^{m_2})$ and $O((L-t)^{m_1})$; boundedness of $F$ makes the boundary terms in the first-order distributional derivative identity tend to zero. Thus the extension is globally $W^{1,4}$.
In the Euler equation, the integrated radial flux, apart from a fixed angular volume factor, is $\frac{P}{k}\partial_s$ of target norm $|P|$, which tends to zero by \eqref{eq:zero-end-flux}. The target curvature term in the extrinsic equation is bounded by $C|dF|^4$. Integration by parts and dominated convergence give the weak equation for arbitrary target variations across both strata.
\end{proof}

\begin{remark}[A concrete pair of nonellipsoidal smooth geometries]\label{new:example-lapses}
Take round sphere factors, $a=b=c=d=1$, and
\[
 h(t)^2=1+\delta\sin^2(2t),\quad \delta>-1,
 \qquad k(s)^2=\frac1{1+\varepsilon\sin^2s\cos^2s},\quad 0\le\varepsilon<1.
\]
These coefficients are even at the endpoints and equal to one there; in radial Cartesian coordinates the metrics extend smoothly over the collapsed sphere factors. For example, at zero $r=\sin t$ makes the radial coefficient $\frac{h(t)^2}{1-r^2}$ a smooth function of $r^2$ equal to one at $r=0$; the other coefficients are also smooth functions of $r^2$. Put $p(v)=v(1-v)$. Then the quantity inside the power $\frac{2}{3}$ in \eqref{new:eq:beta-k} is $r(v)=p(v)+\varepsilon p(v)^2$, with
\[
 r''(v)=-2+\varepsilon(2-12v+12v^2)<0.
\]
Thus $\beta_k=r^{\frac{2}{3}}$ is concave. Theorem~\ref{new:lapses} gives the continuous globally weakly symphonic minimizing join, for example with identity angular maps. Nonzero parameters change the radial geometries beyond the induced ellipsoid family. This example does not assert that the weak minimizer is smooth at the endpoints.
\end{remark}
For arbitrary closed angular factors the theorem still proves all the scalar and regular-region assertions. Such factors generally produce singular metric joins. Proposition~\ref{g:prop:link} explains why that construction cannot be described as a smooth closed-manifold theorem without the completion hypothesis.

\begin{corollary}[Strict stability under radial variations]\label{cor:radial-stability}
For every $\xi\in C_c^\infty(I)$,
\begin{equation}\label{eq:strict-radial-stability}
 d^2J(\widehat\varphi)[\xi,\xi]
 \ge\frac12\int_I\rho(A+B)\sin^2(2\widehat\varphi)\xi^2\,dt.
\end{equation}
The right side is positive when $\xi\not\equiv0$.
\end{corollary}
\begin{proof}
For sufficiently small $\varepsilon$, $\widehat\varphi+\varepsilon\xi\in X_0$.
Apply \eqref{eq:main-energy-gap}, use $dJ(\widehat\varphi)[\xi]=0$, and compare coefficients of $\varepsilon^2$.
\end{proof}
In this inequality the angular maps are fixed.

\section{Ellipsoidal join corollaries}\label{sec:join-equations}
For this specialization, $a,b,c,d>0$ and $m_i,n_i\ge1$. Write
\begin{equation}\label{eq:concise:002}
L=\frac\pi2,\quad I=(0,L),\quad
h(t)=\sqrt{a^2\sin^2t+b^2\cos^2t},\quad
k(s)=\sqrt{c^2\sin^2s+d^2\cos^2s},
\end{equation}
\begin{equation}\label{eq:concise:003}
\rho(t)=h(t)\cos^{m_1}t\sin^{m_2}t.
\end{equation}
The source ellipsoid is
\begin{equation}\label{eq:concise:004}
Q^{m_1+m_2+1}(a,b)
=\{(X,Y)\in\R^{m_1+1}\times\R^{m_2+1}:\frac{|X|^2}{a^2}+\frac{|Y|^2}{b^2}=1\}.
\end{equation}
Its regular region is parametrized by $(t,x,y)\mapsto(a\cos t\,x,b\sin t\,y)$. A prime on $k$ denotes differentiation in its argument $s$.
For $f_i:S^{m_i}\to S^{n_i}$, set
\begin{equation}\label{eq:concise:005}
Q_i=\|f_i^*g_{S^{n_i}}\|^2,\quad V_i=\Vol(S^{m_i}),\quad
C_0=a^{m_1}b^{m_2}V_1V_2,
\end{equation}
\begin{equation}\label{eq:concise:006}
a_1=\frac{c^4F(f_1)}{a^4V_1},\qquad
a_2=\frac{d^4F(f_2)}{b^4V_2}.
\end{equation}
The coefficients $a_1,a_2$ are distinct from the semiaxes $a,b$.

Write $M=Q^{m_1+m_2+1}(a,b)$ and $N=Q^{n_1+n_2+1}(c,d)$ for the join source and target. The regular region is $M_{\mathrm{reg}}=\{0<t<L\}$.
Use the functional and space in \eqref{eq:J}--\eqref{eq:X} with these coefficients.
For clarity, the weak equation used here is the extrinsic Euler equation for the ellipsoid. Let
$\mathsf D=\operatorname{diag}(c^{-2}I_{n_1+1},d^{-2}I_{n_2+1})$.
A map $\Phi\in W^{1,4}(M,N)$ is called weakly symphonic if, with $S=\Phi^*g_N$,
\begin{equation}\label{eq:global-weak-definition}
\begin{aligned}
 \int_M S^{ij}\langle\partial_i\Phi,\partial_j\eta\rangle\,dv_g
 &=\int_M\Lambda\langle\mathsf D\Phi,\eta\rangle\,dv_g
 \quad\text{for every }\eta\in C^\infty(M,\R^{n_1+n_2+2}),\\
 \Lambda&=\frac{S^{ij}\langle\mathsf D\partial_i\Phi,\partial_j\Phi\rangle}
 {|\mathsf D\Phi|^2}.
\end{aligned}
\end{equation}
Here both indices of $S^{ij}$ are raised with $g$; $S^{ij}$ is not an inverse tensor. This is the weak target Euler equation, without an additional assertion of stationarity under domain variations.

\subsection{Geometry and extrinsic formulation}
Differentiating the ellipsoidal parametrization gives the metric and volume form
\begin{equation}
g=a^2\cos^2t\,g_{S^{m_1}}+b^2\sin^2t\,g_{S^{m_2}}+h(t)^2dt^2,
\label{eq:metric}
\end{equation}
\begin{equation}
dv_g=a^{m_1}b^{m_2}\rho(t)\dd\,dv_{S^{m_1}}\,dv_{S^{m_2}}.
\label{eq:volume}
\end{equation}
With our normalization, the standard first variation is
\begin{equation}
\delta E(f)[V]=-2\int_M\langle V,\Div_g\sigma_f\rangle\,dv_g
\label{eq:first}
\end{equation}
for compactly supported $V$; see \cite{KN,NT}.

For scalar functions, define
\begin{equation}\label{eq:operator}
\Delta_g^S v=\frac1{\sqrt{\det g}}\partial_i
 (\sqrt{\det g}\,S^{ij}\partial_jv).
\end{equation}
Suppose $N$ is isometrically embedded in Euclidean space, $S=U^*g_N$, and $D$ is the Euclidean connection.
For fields tangent to $N$, $\nabla^N=(D)^\top$. The Gauss formula therefore gives
\cite[Chapter 8]{Lee}
\begin{equation}\label{eq:extrinsic-projection}
\operatorname{div}_g\sigma_U=(\Delta_g^S U)^\top,
\qquad U\text{ symphonic}\ \Longleftrightarrow\ (\Delta_g^S U)^\top=0.
\end{equation}
Here \eqref{eq:operator} acts componentwise on Euclidean coordinates, and $\top$ is orthogonal projection onto the target tangent space.
For example, a unit sphere map $f$ satisfies
\begin{equation}\label{eq:sphere-extrinsic-check}
\tfrac12\Delta_g^H|f|^2
 =\langle\Delta_g^Hf,f\rangle+H^{ij}\langle f_i,f_j\rangle=0,
\qquad H=f^*g_{S^r}.
\end{equation}
Since $H^{ij}\langle f_i,f_j\rangle=|H|_g^2$,
its tangential projection is $\Delta_g^Hf+|H|_g^2f$.
Thus
\begin{equation}\label{eq:sphere}
f\text{ symphonic}\quad\Longleftrightarrow\quad
\Delta_g^{f^*g_{S^n}}f+\|f^*g_{S^n}\|_g^2f=0.
\end{equation}

For an ellipsoid $N=\{U:U^T\mathsf DU=1\}$, a nonunit normal is $\mathsf DU$.
Applying the product rule for $\Delta_g^S$ to the constraint gives
$\langle\Delta_g^S U,\mathsf DU\rangle
 =-S^{ij}\langle \mathsf DU_i,U_j\rangle$.
Thus the extrinsic equation is
\begin{equation}\label{eq:ellipsoid}
\begin{aligned}
\Delta_g^S U+\Lambda \mathsf DU&=0,\\
\Lambda&=\frac{S^{ij}\langle \mathsf DU_i,U_j\rangle}{|\mathsf DU|^2}.
\end{aligned}
\end{equation}

\subsection{Two source spheres and two target blocks}
This subsection additionally specializes the source to two sphere factors and takes $f_i=f_i(x_i)$, so $r=m=2$ here. It is a special case of Corollaries~\ref{blocks:cor:fixed} and \ref{multi:cor:ellipsoid}; the general two-block target imposes no such source restriction. The following formulas record the coefficients needed for the scalar existence and endpoint results.
Set
$\Phi=(c\cos\varphi\,f_1,d\sin\varphi\,f_2)$,$H_i=f_i^*g_{S^{n_i}}$.
Since $\langle f_i,df_i\rangle=0$,
\begin{equation}\label{eq:concise:014}
S=c^2\cos^2\varphi\,H_1+d^2\sin^2\varphi\,H_2
+k(\varphi)^2(\varphi')^2dt^2,
\end{equation}
\begin{equation}
\|S\|_g^2=\frac{k(\varphi)^4}{h^4}(\varphi')^4
+\frac{c^4\cos^4\varphi}{a^4\cos^4t}Q_1(x)
+\frac{d^4\sin^4\varphi}{b^4\sin^4t}Q_2(y).
\label{eq:join-norm}
\end{equation}
\begin{corollary}\label{lem:join-equation}
If $f_1,f_2$ are symphonic and $\varphi\in C^2(I)$, the join satisfies the symphonic equation on the regular region if and only if
\begin{equation}
\begin{split}
0={}&\frac1\rho\left(\rho\frac{k^4}{h^4}(\varphi')^3\right)'
-\frac{k^3k'}{h^4}(\varphi')^4\\
&+\frac{c^4Q_1(x)}{a^4\cos^4t}\cos^3\varphi\sin\varphi
-\frac{d^4Q_2(y)}{b^4\sin^4t}\sin^3\varphi\cos\varphi
\end{split}
\label{eq:join-pointwise}
\end{equation}
pointwise. Here and below, $k,k',k''$ with omitted arguments are evaluated at $\varphi(t)$.
\end{corollary}
\begin{proof}
Apply Corollary~\ref{g:thm:reduction} with the angular scales $c|\cos s|$, $d|\sin s|$, $a\cos t$, $b\sin t$. On each target chart where $\sin s\cos s\ne0$, absorb the fixed signs into the sphere maps. This gives \eqref{eq:join-pointwise}. At a degenerate latitude the same ambient identity extends by continuity; on any interval where $\varphi$ is constant at that latitude, the map equation is the remaining factor equation. Thus the equivalence holds for the stated $C^2$ functions $\varphi$.
\end{proof}
\begin{corollary}\label{cor:constant}
Under Corollary \ref{lem:join-equation}, if the join is symphonic and $0<\varphi(t_0)<L$ for some $t_0\in I$, then $Q_1$ and $Q_2$ are individually constant.
\end{corollary}
\begin{proof}
At an interior target latitude both angular warps have nonzero derivative. Apply Corollary~\ref{g:cor:constant}.
\end{proof}

For constant densities, the equation reduces to
\begin{equation}
\frac1\rho\left(\rho\frac{k^4}{h^4}(\varphi')^3\right)'
-\frac{k^3k'}{h^4}(\varphi')^4
+a_1\frac{\cos^3\varphi\sin\varphi}{\cos^4t}
-a_2\frac{\sin^3\varphi\cos\varphi}{\sin^4t}=0.
\label{eq:join-radial}
\end{equation}
Equivalently, its expanded form is
\begin{equation}
\begin{split}
0={}&\frac{k^4}{h^4}\bigl((\varphi')^3\bigr)'
+\frac{k^4}{h^4}\left(-m_1\tan t+m_2\cot t-3\frac{h'}h\right)(\varphi')^3
\\&+\frac{3k^3k'}{h^4}(\varphi')^4
+a_1\frac{\cos^3\varphi\sin\varphi}{\cos^4t}
-a_2\frac{\sin^3\varphi\cos\varphi}{\sin^4t}.
\end{split}
\label{eq:radial-expanded}
\end{equation}

\begin{lemma}\label{lem:energy}
Even for nonconstant $Q_i$, the averaged coefficients in \eqref{eq:concise:006} give
\begin{equation}
F(\phi)=4C_0J(\varphi),\qquad E(\phi)=2C_0J(\varphi).
\label{eq:energy-reduction}
\end{equation}
\end{lemma}
\begin{proof}
Integrate \eqref{eq:join-norm} using \eqref{eq:volume} and $\int Q_i=F(f_i)$. Nonnegativity makes the identity valid also for infinite energy.
\end{proof}
This identity concerns angular averages only; it does not permit replacing the pointwise densities in \eqref{eq:join-pointwise} by their averages.

Differentiating \eqref{eq:J} on compact interior intervals gives the first variation with respect to $\varphi$. Under the factor hypotheses of Corollary \ref{thm:join}, its Euler equation is \eqref{eq:join-radial}, so Corollary \ref{lem:join-equation} yields the full map equation. Without constant densities, variations of $\varphi$ yield only the angular average of the radial equation.
\begin{corollary}[Minimizing join]\label{thm:join}
Suppose $f_i$ are smooth and symphonic, with $Q_i\equiv \frac{F(f_i)}{V_i}>0$. The functional $J$ on $X_0$ has a unique finite-energy minimizer. Its representative $\widehat\varphi$ satisfies
\begin{equation}\label{eq:concise:008}
\widehat\varphi\in C([0,L])\cap C^\infty(I),\qquad
\widehat\varphi(0)=0,\quad\widehat\varphi(L)=L,
\qquad \widehat\varphi'>0\text{ on }I.
\end{equation}
The corresponding join
\begin{equation}\label{eq:concise:009}
\phi(a\cos t\,x,b\sin t\,y)
=\bigl(c\cos\widehat\varphi(t)f_1(x),d\sin\widehat\varphi(t)f_2(y)\bigr)
\end{equation}
belongs to $W^{1,4}(M,N)\cap C^0(M,N)$, is weakly symphonic on $M$, is smooth and symphonic on $M_{\mathrm{reg}}$, and satisfies $E(\phi)=2C_0J(\widehat\varphi)$. If it extends smoothly across both endpoints, that extension is symphonic everywhere.
\end{corollary}
\begin{proof}
The coefficients in \eqref{eq:concise:002} are positive, smooth and even at both endpoints. Put
\[
 r(v)=\frac{v(1-v)}{d^2+(c^2-d^2)v}.
\]
Then $r''=-\frac{2c^2d^2}{[d^2+(c^2-d^2)v]^3}<0$, so $\beta_k=r^{\frac{2}{3}}$ is concave. The angular hypotheses are precisely those of Theorem~\ref{new:lapses}. The completions are the smooth induced ellipsoids, and \eqref{eq:energy-reduction} fixes the energy normalization. Apply that theorem; a smooth extension satisfies the equation everywhere by continuity.
\end{proof}
\begin{corollary}[Quantitative isolation on ellipsoids]\label{thm:energy-gap}
Under Corollary~\ref{thm:join}, the estimate \eqref{eq:main-energy-gap} holds for every finite-energy competitor in $X_0$. Every minimizing sequence converges strongly in $X$ and locally uniformly, and \eqref{eq:strict-radial-stability} holds.
\end{corollary}
\begin{proof}
These are the gap and convergence assertions of Theorem~\ref{new:lapses} and Corollary~\ref{cor:radial-stability}, under the same specialization.
\end{proof}
\begin{remark}\label{rem:constant-functions}
The constant functions $\varphi\equiv0$ and $\varphi\equiv L$ have energies
\begin{equation}\label{eq:concise:037}
J(0)=\frac{a_1}4\int_Ih\cos^{m_1-4}t\sin^{m_2}t\dd,\qquad
J(L)=\frac{a_2}4\int_Ih\cos^{m_1}t\sin^{m_2-4}t\dd.
\end{equation}
The first is finite exactly when $m_1>3$, the second exactly when $m_2>3$. The obstacle variation in the proof of Theorem~\ref{new:lapses} shows neither can be a finite energy minimizer. These integrability thresholds do not restrict Corollary \ref{thm:join}.
\end{remark}

\section{Smooth endpoint obstructions}\label{sec:endpoints}
Near $t=0$, write $z=r\omega$ in the collapsing source directions. Smoothness of a join requires the two expressions
\begin{equation}\label{eq:end-smooth}
c\cos\varphi(\arcsin\left(\frac{r}{b}\right))f_1(x),\qquad
d\sin\varphi(\arcsin\left(\frac{r}{b}\right))f_2\left(\frac{z}{r}\right)
\end{equation}
to be smooth in Cartesian coordinates at $z=0$. The endpoint values alone do not imply this condition.
\subsection{Two endpoint estimates for the general rigidity theorem}

\begin{lemma}[Exclusion of vanishing beyond second order]\label{lem:no-flat-endpoint}
Let $a_1,a_2>0$. Suppose $\varphi\in C([0,L])\cap C^2(I)$ satisfies \eqref{eq:join-radial}, takes values in $[0,L]$, and has endpoint values $0,L$.
If $\varphi$ has a $C^2$ one-sided extension at zero, then $\varphi(t)=o(t^2)$ is impossible.
The corresponding assertion holds for $L-\varphi(L-t)$ at the other endpoint.
\end{lemma}
\begin{proof}
Suppose $\varphi=o(t^2)$. The $C^2$ extension gives $\varphi'=o(t)$.
With $u=K(\varphi)$ and $P=q(u')^3$, we have $P(0+)=0$.
The moving minimum satisfies
\begin{equation}\label{eq:general-moving-well}
 \alpha_*(t)=\arctan\left(\sqrt{\frac{a_1}{a_2}}\tan^2t\right)
 \sim\sqrt{\frac{a_1}{a_2}}\,t^2.
\end{equation}
Thus $0\le\varphi<\alpha_*$ near zero, and \eqref{eq:flux} gives $P'\le0$.
It follows that $P\le0$ and $u'\le0$. Since $u\ge0$ and $u(0)=0$, we obtain $u=0$ near zero.

The zero set of $u$ is relatively closed in $I$.
If $u(t_0)=0$, then $u'(t_0)=P(t_0)=0$ and $u<u_*$ on a neighborhood of $t_0$.
There $P'\le0$, hence $u'\ge0$ to the left and $u'\le0$ to the right.
Nonnegativity forces $u=0$ on the whole neighborhood.
The zero set is therefore also relatively open, so $u\equiv0$ on $I$, contrary to the right endpoint value.
The right endpoint assertion follows by exchanging the two factors.
\end{proof}

Recall that a degree-$\ell$ spherical eigenmap is a map $f:S^m\to S^n$ whose components are restrictions of homogeneous harmonic polynomials of degree $\ell$; the polar-coordinate formula for the Euclidean Laplacian gives the eigenvalue below.
Its components satisfy
\begin{equation}\label{eq:eigenmap-eigenvalue}
 \Delta_{S^m}f=-\ell(\ell+m-1)f,
 \qquad |df|^2=\ell(\ell+m-1).
\end{equation}
The second identity follows by applying the Laplacian to $|f|^2=1$.

\begin{lemma}[Density bound for a quadratic eigenmap]\label{lem:quadratic-density-bound}
If $f:S^m\to S^n$ is a degree-two spherical eigenmap, then
\begin{equation}\label{eq:quadratic-density-bound}
 \|f^*g_{S^n}\|^2\le8(m+1)
\end{equation}
pointwise.
\end{lemma}
\begin{proof}
Write the homogeneous polynomial map as $F(y)=\mathcal B(y,y)$, where $\mathcal B$ is a symmetric vector-valued bilinear form.
Homogeneity and $|f|=1$ imply $|F(z)|^2=|z|^4$ for every $z\in\R^{m+1}$.
Fix orthonormal vectors $y,\xi$ and put
$A=\mathcal B(y,y)$, $D=\mathcal B(\xi,\xi)$, and $C=\mathcal B(y,\xi)$.
Comparing the $s^2t^2$ coefficients in
$|s^2A+2stC+t^2D|^2=(s^2+t^2)^2$ gives
\begin{equation}\label{eq:quadratic-bilinear-identity}
 2\langle A,D\rangle+4|C|^2=2.
\end{equation}
Since $|A|=|D|=1$, we have $|C|\le1$ and $|df_y(\xi)|=2|C|\le2$.
Let $\lambda_1,\ldots,\lambda_m$ be the eigenvalues of the pullback metric relative to the unit sphere metric.
Then $0\le\lambda_j\le4$ and, by \eqref{eq:eigenmap-eigenvalue}, $\sum_j\lambda_j=2(m+1)$.
Consequently $\sum_j\lambda_j^2\le4\sum_j\lambda_j=8(m+1)$.
\end{proof}

\subsection{Angular linearity and its specializations}
\begin{theorem}[Every smooth connecting join has linear angular factors]\label{new:linear}
Let $f_i:S^{m_i}\to S^{n_i}$ be smooth nonconstant maps, with $m_i,n_i\ge1$. Let $\varphi\in C([0,L])\cap C^2(I)$ take values in $[0,L]$ and connect $0$ to $L$. Suppose the associated ellipsoidal join is symphonic on the regular region and extends smoothly to the entire source ellipsoid. Then each $f_i$ is the restriction of a linear isometric embedding. This holds for every $a,b,c,d>0$.
\end{theorem}
\begin{proof}
There is a $t_0$ with $0<\varphi(t_0)<L$. The full angular equation in Corollary~\ref{g:thm:reduction}, applied on a neighborhood of that slice with the ellipsoidal warps, implies that each $f_i$ is symphonic. Subtraction of the full radial equation at two points of one factor, with the other variable fixed, gives constant densities by Corollary~\ref{cor:constant}. Nonconstancy of $f_i$ implies $Q_i>0$. The radial equation \eqref{eq:join-radial} holds on the open set $0<\varphi<L$. On the interior of each pole level set $\varphi$ is constant and each term in that equation vanishes. Continuity gives the equation at the remaining points as well.

Consider the second factor, and put $m=m_2$. Fix the first angular variable. In normal Cartesian source coordinates $z=r\omega$, the second Euclidean target block is
\[
 G(z)=\eta(r)f_2(\omega),\qquad \eta(r)=d\sin\varphi(\arcsin\left(\frac{r}{b}\right)).
\]
It is smooth at $z=0$. Along one fixed ray, taking the inner product with the unit vector $f_2(\omega_0)$ recovers $\eta(r)$, so $\eta$ and then $\varphi$ have smooth one-sided extensions. By Lemma~\ref{lem:no-flat-endpoint}, $\varphi=o(t^2)$ is impossible. Its first nonzero Taylor order is therefore $\ell=1$ or $2$, with positive coefficient. Write $\eta(r)=Cr^\ell+o(r^\ell)$, $C>0$.

Since $|G(z)|=|\eta(|z|)|=O(|z|^\ell)$, all Taylor polynomials of degree below $\ell$ vanish. Its homogeneous term $P_\ell$ of degree $\ell$ satisfies
\[
 P_\ell(\omega)=\lim_{r\downarrow0}r^{-\ell}G(r\omega)=Cf_2(\omega).
\]
Thus the polynomial character of the factor follows from smoothness of the whole map. If $\ell=1$, $f_2(\omega)=A\omega$ and $|A\omega|=1$ for every unit $\omega$, giving $A^TA=I$.

If $\ell=2$, integrate the sphere-valued symphonic equation, componentwise,
\[
 \operatorname{div}_{S^m}(S^{ij}\partial_j f_2)+Q_2f_2=0.
\]
The divergence integrates to zero and $Q_2$ is a positive constant, so $\int_{S^m}f_2=0$. A homogeneous quadratic polynomial $\mathcal B(z,z)$ has sphere mean proportional to $\tr\mathcal B$; consequently every component of $P_2$ is harmonic. The map $f_2$ is therefore a degree-two spherical eigenmap. Smoothness along a signed ray and $f_2(-\omega)=f_2(\omega)$ make $\varphi$ even across $t=0$, so
$\varphi(t)=Ct^2+O(t^4)$ with $C>0$. Substitution into \eqref{eq:join-radial} gives
\begin{equation}\label{eq:quadratic-general-coefficient}
 0=\left[a_1C+\frac{d^4}{b^4}\{8(m+3)-Q_2\}C^3\right]t^2+O(t^4).
\end{equation}
The divergence contributes $\frac{8(m+3)d^4C^3}{b^4}$, the $k'$ term is higher order, and the two angular terms give the other displayed coefficients. Lemma~\ref{lem:quadratic-density-bound} yields $8(m+3)-Q_2\ge16$, so the leading coefficient is strictly positive. This contradiction leaves only $\ell=1$. Reflection about the other endpoint proves the same assertion for $f_1$.
\end{proof}
Proposition~\ref{prop:homothety} gives a smooth example when the source and target semiaxes are proportional.

\begin{corollary}[Obstruction for nonlinear eigenmap factors]\label{thm:eigen-obstruction}
Suppose $f_i:S^{m_i}\to S^{n_i}$ are smooth symphonic maps with positive constant angular densities. If either factor is the restriction of a vector of homogeneous harmonic polynomials of a common degree $\ell\ge2$, there is no globally smooth symphonic join of these factors with
$0\le\varphi\le L$, $\varphi(0)=0$, and $\varphi(L)=L$.
The conclusion holds for every $a,b,c,d>0$.
\end{corollary}
\begin{proof}
Theorem~\ref{new:linear} forces every nonconstant factor to be linear isometric. A spherical harmonic of degree at least two cannot be a nonzero degree-one spherical harmonic. This excludes the stated factor.
\end{proof}
\begin{corollary}\label{cor:degree-one-necessary}
If both angular factors are spherical eigenmaps of positive degree satisfying the hypotheses of Corollary \ref{thm:join}, a globally smooth symphonic join can exist only when both degrees are one.
Each factor is then the restriction of a linear isometric embedding.
\end{corollary}
\begin{proof}
Corollary \ref{thm:eigen-obstruction} excludes every degree at least two.
A degree-one factor has the form $f(y)=Ay$; the identity $|Ay|=1$ on $S^m$ implies $A^TA=I$.
\end{proof}

\subsection{Examples and necessity of the hypotheses}
\begin{corollary}[The Veronese family]\label{prop:veronese-general}
Let $m=m_2\ge1$, and take the quadratic Veronese map as the second factor:
\begin{equation}\label{eq:general-veronese}
f_2(y)=\sqrt{\frac{m+1}{m}}\left(yy^T-\frac{I_{m+1}}{m+1}\right)
 \in S^{\frac{(m+1)(m+2)}{2}-2},\qquad y\in S^m,
\end{equation}
The space of traceless symmetric matrices carries the Frobenius inner product. Assume the first factor $f_1$ is smooth and symphonic with positive constant density $Q_1$.
For arbitrary $a,b,c,d>0$, no globally smooth symphonic join exists with $0\le\varphi\le L$, $\varphi(0)=0$, and $\varphi(L)=L$. Interchanging the factors gives the right endpoint statement.
\end{corollary}
\begin{proof}
Direct calculation gives $|f_2|=1$. For a unit tangent vector $\xi\perp y$,
$df_2(\xi)=\sqrt{\frac{m+1}{m}}(\xi y^T+y\xi^T)$, hence
\begin{equation}\label{eq:veronese-density-general}
 f_2^*g=\frac{2(m+1)}m g_{S^m},\qquad Q_2=\frac{4(m+1)^2}m,
 \qquad \Delta f_2=-2(m+1)f_2.
\end{equation}
Thus $f_2$ is harmonic; because its pullback tensor is a constant multiple of the source metric, it is also symphonic.
The components are degree-two spherical harmonics. Corollary \ref{thm:eigen-obstruction} therefore applies.
\end{proof}
\begin{proposition}\label{prop:nonconstant}
Symphonicity of the factors alone does not imply existence of a symphonic join connecting $0$ to $L$.
\end{proposition}
\begin{proof}
Take $m_1=n_1=7$, $m_2=n_2=4$, $a=b=c=1$, $d=\frac{1}{2}$, and $f_1=\id_{S^7}$. In stereographic coordinates on $S^4$, let $f_2(z)=2z$. This smooth conformal diffeomorphism has
\begin{equation}\label{eq:concise:039}
f_2^*g_{S^4}=\lambda(z)^2g_{S^4},\qquad
\lambda(z)=\frac{2(1+|z|^2)}{1+4|z|^2}.
\end{equation}
If $\widehat g=e^{2w}g$ on a four-dimensional source, then $|v^*G|_{\widehat g}^2=e^{-4w}|v^*G|_g^2$ and $dv_{\widehat g}=e^{4w}dv_g$. Change of variables therefore gives $E(v\circ\psi)=E(v)$ for every conformal diffeomorphism $\psi$ of the source. Arbitrary variations correspond through $\psi^{-1}$, so conformal precomposition preserves criticality. As the identity is symphonic, so is $f_2$, and
\begin{equation}\label{eq:concise:040}
Q_1=7,\quad Q_2=4\lambda^4\text{ nonconstant},\quad
F(f_2)=4V_2,\quad a_1=7,\quad a_2=\frac14.
\end{equation}
However, any continuous $\varphi$ connecting $0$ to $L$ takes an interior value, and Corollary \ref{cor:constant} rules out a symphonic join.
\end{proof}

\begin{proposition}\label{prop:homothety}
Let $n_i=m_i$, $f_i=\id_{S^{m_i}}$, $c=\lambda a$, and $d=\lambda b$, with $\lambda>0$. Then $\varphi(t)=t$ yields the globally smooth symphonic map $\phi(z)=\lambda z$.
\end{proposition}
\begin{proof}
We have $\phi^*g_N=\lambda^2g$. Constant rescaling preserves the Levi--Civita connection \cite[Chapter 5]{Lee}, so $\nabla d\phi=0$ and $\Div\sigma_\phi=\lambda^2\tau(\phi)=0$. The ambient linear formula guarantees endpoint smoothness.
\end{proof}

\section{Hopf constructions: general warps and the full equation}\label{hc:sec:setup}
Fix an integer $r\ge2$. The Hopf construction uses one angular map and a radial function $u$ connecting the poles. We use the normalization $\mathcal E=\frac{E}{2}$, with $E$ as in \eqref{eq:energy}; this positive factor changes no critical points. The variable $u$ is target arclength, with poles at $0,U$, whereas the later ellipsoidal specialization uses a latitude variable with poles at $-\frac{\pi}{2},\frac{\pi}{2}$. The normalized $J$ below also removes the fixed angular volume factor.\par\medskip For the local equation let $(P_i,g_i)$ be arbitrary connected Riemannian manifolds, $m_i=\dim P_i\ge1$, for $1\le i\le r$, and let $(N,\gamma)$ be an arbitrary connected Riemannian manifold. No compactness or curvature hypothesis is needed for that equation. For the global variational results, the normalized energy formula and the statements on metric completion below, assume these angular manifolds are closed. On the regular source and target consider
\begin{align}
 M_{\rm reg}&=(0,T)\times\prod_{i=1}^rP_i,
 &g&=h(t)^2dt^2+\sum_{i=1}^ra_i(t)^2g_i,\label{hc:eq:source}\\
 Y_{\rm reg}&=(0,U)\times N,
 &G&=du^2+b(u)^2\gamma.\label{hc:eq:target}
\end{align}
Here $T,U>0$, the source functions are smooth and positive in the interior, and $b\in C^\infty([0,U])$ satisfies
\begin{equation}\label{hc:eq:poles}
 b(0)=b(U)=0,\qquad b>0\text{ on }(0,U),\qquad b'(0)>0,\quad b'(U)<0.
\end{equation}
Writing the target radial coefficient as one loses no generality for a metric $k(s)^2ds^2+B(s)^2\gamma$: take $u=\int k(s)ds$ and $b(u)=B(s(u))$, provided this arclength interval is finite.

Write $P=\prod_{i=1}^rP_i$ and $z=(x_1,\ldots,x_r)$. The single map $f:P\to N$ may depend on every factor and need not split into separate factor maps. In particular, no spherical or ellipsoidal assumption is imposed on the $P_i$ here. The Hopf construction with radial component $u(t)$ is
\begin{equation}\label{hc:eq:ansatz}
 H_f(t,x_1,\ldots,x_r)=(u(t),f(x_1,\ldots,x_r)),\qquad u(0)=0,\quad u(T)=U.
\end{equation}
The pole limits define a continuous map to the suspension of $N$ after prescribed source endpoint collapses. More precisely, at each endpoint one may collapse a selected subset of angular factors and retain the other factors; the map sends the whole endpoint stratum to the corresponding target pole, so it respects these identifications. For $r=2$, collapsing $P_2$ at zero and $P_1$ at $T$ recovers $P_1*P_2$. For $r>2$ we retain one interval variable: this is not the full $r$-fold topological join, which has $r-1$ joining parameters. The quotient agrees with the metric completion only under the corresponding finite-distance collapse conditions. The completed spaces need not be smooth manifolds for arbitrary links. This is the same separated construction as the sphere and ellipsoid Hopf ansatz, with the target written in polar coordinates; see \cite{Cao,DFL} for the symphonic and harmonic settings, respectively.

We use $\mathcal E(F)=\tfrac14\int|F^*G|_g^2dv_g$ and
\[
 \sigma_F(X)=\sum_jG(dF(X),dF(e_j))dF(e_j),\qquad \cT(F)=\divg_g\sigma_F.
\]
Set $g_t=\sum_{i=1}^ra_i^2g_i$, $\rho=h\prod_{i=1}^ra_i^{m_i}$, and
\begin{equation}\label{hc:eq:qt}
 Q_t(z)=|f^*\gamma|_{g_t}^2
 =\sum_{i=1}^r\frac{Q_i(z)}{a_i^4}
       +2\sum_{1\le i<j\le r}\frac{Q_{ij}(z)}{a_i^2a_j^2}.
\end{equation}
In fixed $g_i$-orthonormal frames $e_{i\alpha}$, these densities are
\[
 Q_i=\sum_{\alpha,\beta}\gamma(df(e_{i\alpha}),df(e_{i\beta}))^2,
 \qquad
 Q_{ij}=\sum_{\alpha,\beta}\gamma(df(e_{i\alpha}),df(e_{j\beta}))^2.
\]
Thus all diagonal blocks and all pairwise mixed blocks are retained. For $r=3$, the mixed terms are $\frac{2Q_{12}}{a_1^2a_2^2}$, $\frac{2Q_{13}}{a_1^2a_3^2}$ and $\frac{2Q_{23}}{a_2^2a_3^2}$; none may be discarded without justification.

\begin{corollary}[Full Hopf-construction criterion]\label{hc:prop:full}
For $u\in C^2((0,T),(0,U))$, the map \eqref{hc:eq:ansatz} is symphonic on the regular region if and only if
\begin{align}
 \cT_{g_t,\gamma}(f)&=0\quad\text{for every }t,\label{hc:eq:angular}\\
 \frac1\rho\left(\frac\rho{h^4}u'^3\right)'
             -b(u)^3b'(u)Q_t(z)&=0\quad\text{at every }(t,z).\label{hc:eq:radial}
\end{align}
In particular, at each $t$ for which $b'(u(t))\ne0$, the radial equation forces $Q_t$ to be independent of $z$.
\end{corollary}
\begin{proof}
This is the $m=1$ specialization of Corollary~\ref{coupled:cor:warped}, with arbitrary source-factor count $r$ and $\partial_tf_1=0$. Equivalently, in Corollary~\ref{new:master}, take $k=1$, $H_u=b(u)^2\gamma$ and $g_t=\sum a_i(t)^2g_i$. Then $Q=b^4Q_t$, the angular tension is $b^2\mathcal T_{g_t,\gamma}(f)$, and $\frac{Q_u}{4}=b^3b'Q_t$. Positivity of $b$ gives the equivalence; subtraction at two angular points gives the final assertion.
\end{proof}

For the existence theorems we impose the following explicit compatibility hypothesis:
\begin{equation}\label{hc:eq:compatibility}
 \cT_{g_t,\gamma}(f)=0,\qquad Q_t(z)=q(t)>0\quad\text{for every }t,z.
\end{equation}
It can be checked without solving the radial equation. Requiring the separate slices of $f$ to be symphonic is not, on its own, a substitute for \eqref{hc:eq:compatibility}. Likewise, replacing $Q_t$ by its angular average does not give the full map equation. Constancy of every $Q_i$ and $Q_{ij}$ is a useful sufficient condition for the density part, but the angular equation still has to be verified.

With
\begin{equation}\label{hc:eq:weights}
 A(t)=\frac{\rho(t)}{h(t)^4},\qquad C(t)=\rho(t)q(t),\qquad V(u)=b(u)^4,
\end{equation}
the reduced energy is
\begin{equation}\label{hc:eq:energy}
 J(u)=\frac14\int_0^T\bigl(A|u'|^4+CV(u)\bigr)dt,
 \qquad\mathcal E(H_f)=\left(\prod_{i=1}^r\Vol(P_i)\right)J(u).
\end{equation}
The admissible class $\mathcal A$ consists of continuous functions on $[0,T]$, locally $W^{1,4}$ in the interior, with $0\le u\le U$, the two stated traces, and finite $J$. It is nonempty: take a smooth transition on a compact interior interval and the respective pole values near the two endpoints.

\section{Connecting solutions of the radial Hopf equation}\label{hc:sec:exist}
\begin{theorem}[Trace-controlled Hopf minimizers]\label{hc:thm:exist}
For any finite $r\ge2$ as above, assume \eqref{hc:eq:poles}, \eqref{hc:eq:compatibility}, and
\begin{equation}\label{hc:eq:trace}
 \int_0^T A(t)^{-\frac{1}{3}}dt<\infty.
\end{equation}
Suppose that for some $c\in(0,U)$, $b$ is nondecreasing on $[0,c]$ and nonincreasing on $[c,U]$. Then $J$ attains its finite infimum on $\mathcal A$. Every minimizer satisfies
\begin{equation}\label{hc:eq:strict}
 0<u(t)<U,\qquad u'(t)>0\quad(0<t<T),\qquad u\in C^\infty(0,T).
\end{equation}
Consequently \eqref{hc:eq:ansatz} is a nonconstant smooth symphonic Hopf construction on the regular region and has the prescribed continuous pole limits. Uniqueness of the minimizer is not asserted.
\end{theorem}
\begin{proof}
\emph{Existence.} For $0<s<t<T$, weighted H\"older gives
\begin{equation}\label{hc:eq:holder}
 |u(t)-u(s)|\le\left(\int_s^t A|u'|^4\right)^{\frac{1}{4}}
                    \left(\int_s^t A^{-\frac{1}{3}}\right)^{\frac{3}{4}}.
\end{equation}
Thus a minimizing sequence has a common modulus of continuity up to both ends and is uniformly bounded. Pass to a uniformly convergent subsequence on $[0,T]$, and to weakly convergent derivatives in $L^4(A\,dt)$. Testing against compactly supported smooth functions identifies the derivative of the limit. The traces and range persist. Lower semicontinuity controls the kinetic term and Fatou's lemma controls the nonnegative potential. The finite competitor described above proves that the limit is a finite-energy minimizer.

\emph{The equation despite the interval constraint.} Extend $V$ to the real line by even reflection at $0,U$ and period $2U$. Since $V'(0)=V'(U)=0$, this extension is $C^1$. The folding map $R(v)=\min_{j\in\mathbb Z}|v-2jU|$ preserves the potential and the magnitude of the derivative almost everywhere. For a compactly supported variation $\eta$, $R(u+\epsilon\eta)$ is admissible. Hence two-sided variations give
\begin{equation}\label{hc:eq:flux}
 P=Au'^3,\qquad P'=\frac C4V'(u)=Cb(u)^3b'(u)
\end{equation}
in distributions. Its right side is continuous locally, so $P\in C^1$ and $u'=\left(\frac{P}{A}\right)^{\frac{1}{3}}$ gives $u\in C^{1,\frac{1}{3}}_{\rm loc}$.

\emph{Neither pole is reached inside.} If $u(t_0)=0$, nonnegativity gives $u'(t_0)=P(t_0)=0$. On a small compact neighborhood, $|V'(u)|\le K|u|^3$. Set $Z=u^4+|P|^{\frac{4}{3}}$. Equation \eqref{hc:eq:flux}, the local bounds for $A^{-1},C$, and Young's inequality imply
\[
 |Z'|\le K_1|u|^3|P|^{\frac{1}{3}}\le K_2 Z.
\]
Gronwall's inequality in both time directions gives $Z=0$ near $t_0$. The zero set is therefore open and closed in the connected regular interval, forcing $u\equiv0$, contrary to the other trace. The same argument with $U-u$ excludes the upper pole.

\emph{Loop removal proves monotonicity.} Suppose a minimizer decreases somewhere. Choose $s<t$ with $u(s)>u(t)$. Put $d=\max_{[0,t]}u$ and choose $s_0<t$ at which $d$ is attained. Since $0<d<U$, let $r>t$ be the first subsequent time after $t$ at which $u(r)=d$. Let $a=\min_{[s_0,r]}u<d$, attained at $q\in(s_0,r)$, and let $\ell<s_0$ be the last time before $s_0$ with $u(\ell)=a$. This time exists because $u(0)=0<a$. Throughout $[\ell,q]$ and $[s_0,r]$, the relevant values lie in $[a,d]$.

The one-peak hypothesis implies $V(w)\ge\min\{V(a),V(d)\}$ for $w\in[a,d]$. If $V(a)\le V(d)$, replace $u$ on $[\ell,q]$ by the constant $a$. If $V(d)<V(a)$, replace it on $[s_0,r]$ by the constant $d$. In either case the replacement is admissible, the potential does not increase pointwise, and the kinetic term strictly decreases because the removed segment was nonconstant and $A>0$. This contradicts minimality. Thus $u$ is nondecreasing and $P\ge0$.

\emph{Strict increase and smoothness.} If $P(t_0)=0$ and $u(t_0)\le c$, then on $(0,t_0)$ one has $u\le c$, hence $P'\ge0$. Consequently $0\le P(t)\le P(t_0)=0$ there, making $u$ constant on that interval, contrary to $u(0)=0<u(t_0)$. If $u(t_0)\ge c$, use $P'\le0$ on $(t_0,T)$ and the upper trace instead. Therefore $P>0$ everywhere in the interior. The system \eqref{hc:eq:flux} is a smooth ordinary differential system where $P>0$, proving \eqref{hc:eq:strict}. Finally \eqref{hc:eq:compatibility} and Corollary~\ref{hc:prop:full} give the full symphonic map equation.
\end{proof}

\begin{corollary}[The total collapsing dimension controls the trace]\label{hc:cor:subcritical}
Assume $h$ has positive finite endpoint limits. Let $I_0,I_T$ be nonempty subsets of $\{1,\ldots,r\}$. Near zero suppose $a_i(t)\sim\alpha_i t$ for $i\in I_0$, while the other $a_i$ have positive finite limits. Near $T$ suppose $a_i(t)\sim\beta_i(T-t)$ for $i\in I_T$, while the other $a_i$ have positive finite limits. All $\alpha_i,\beta_i$ are positive. Define
\[
 d_0=\sum_{i\in I_0}m_i,\qquad d_T=\sum_{i\in I_T}m_i.
\]
Then the trace condition \eqref{hc:eq:trace} holds if and only if $d_0<3$ and $d_T<3$. Under these inequalities and the remaining hypotheses of Theorem~\ref{hc:thm:exist}, the radial Hopf equation has a strictly increasing solution, smooth in the interior, with the prescribed pole limits. The number of noncollapsing factors is unrestricted.
\end{corollary}
\begin{proof}
At the two endpoints $A$ is asymptotic to positive constants times $t^{d_0}$ and $(T-t)^{d_T}$, respectively. Its inverse cube root is integrable exactly when the corresponding total dimension is less than three. On every compact interior interval it is integrable by positivity and smoothness. For the earlier collapse of two source factors, $I_0=\{2\}$ and $I_T=\{1\}$, so this reduces to $m_1,m_2\in\{1,2\}$. With several simultaneously collapsing factors, their dimensions must be added; checking each factor separately is insufficient.
\end{proof}

\begin{proposition}[Zero flux and the global weak equation]\label{hc:prop:global}
Assume the hypotheses of Theorem~\ref{hc:thm:exist} and the endpoint geometry of Corollary~\ref{hc:cor:subcritical}, with $d_0,d_T<3$. Suppose additionally that
\[
 q(t)\sim q_0t^{-4}\quad(t\downarrow0),\qquad
 q(t)\sim q_T(T-t)^{-4}\quad(t\uparrow T),\qquad q_0,q_T>0.
\]
Then the flux of every minimizer satisfies $P(0+)=P(T-)=0$. If both geometric completions are smooth compact Riemannian manifolds with the indicated collapsed strata, the associated Hopf map extends continuously and its extension is globally $W^{1,4}$ and weakly symphonic under all target variations.
\end{proposition}
\begin{proof}
Near zero write $m=d_0<3$, $A\sim A_0t^m$, $C\sim C_0t^{m-4}$, with positive constants. Since $u\to0$, $b(u)\sim b'(0)u$, $b'(u)\to b'(0)>0$, and $P'=Cb^3b'\ge0$. Thus the positive flux has a finite nonnegative left limit $P_0$. If $P_0>0$, then
\[
 u'\sim\left(\frac{P_0}{A_0}\right)^{\frac{1}{3}}t^{-\frac{m}{3}},\qquad
 u\sim c_0t^{1-\frac{m}{3}},\quad c_0>0.
\]
Substitution in the flux equation gives $P'\sim \frac{c_1}{t}$ with $c_1>0$, contradicting a finite left limit. Therefore $P_0=0$. Reflection gives the other endpoint.

For the weak extension, use a smooth isometric embedding of the compact target. The pullback tensor has nonnegative eigenvalues and $|dH_f|^4\le(\dim M)|H_f^*G|^2$, so the regular-region differential is in $L^4$. The boundary areas tend to zero as $t^{d_0}$ and $(T-t)^{d_T}$, and the map is bounded and continuous with angle-independent pole limits. Integration of the first-order derivative identity on truncated regions therefore identifies its derivatives globally. In the Euler equation, the integrated radial boundary-flux vector is $P\partial_u$, apart from the fixed angular measure. Its norm is $P$, which tends to zero. The target's extrinsic curvature term is bounded by a constant times $|dH_f|^4$, hence integrable. Integration by parts and passage to the limit give the full weak map equation. This argument does not assert endpoint smoothness.
\end{proof}

\section{The radial Hopf equation in the critical case}\label{hc:sec:critical}
The trace integral diverges at a linearly collapsing endpoint of total collapsing dimension three. This does not rule out existence, and the following class handles a critical case exactly. Its target warping need not have one peak.

\begin{theorem}[Exact energy bound and characterization of minimizers]\label{hc:thm:calibration}
For any finite $r\ge2$ as above, assume \eqref{hc:eq:poles} and \eqref{hc:eq:compatibility}. Suppose
\begin{equation}\label{hc:eq:constantD}
 D(t):=A(t)^{\frac{1}{4}}C(t)^{\frac{3}{4}}=\frac{\rho(t)}{h(t)}q(t)^{\frac{3}{4}}=D_0>0
\end{equation}
is constant, and for an interior $t_0$ the function
\[
 Z(t)=\int_{t_0}^t h(\tau)\left(\frac{q(\tau)}3\right)^{\frac{1}{4}}d\tau
\]
has limits $-\infty,+\infty$ at the two endpoints. Define
\[
 \Psi(u)=\int_{\frac{U}{2}}^{u}\frac{d\xi}{b(\xi)}.
\]
Then for every $c\in\R$,
\begin{equation}\label{hc:eq:explicit}
 u_c(t)=\Psi^{-1}(Z(t)+c)
\end{equation}
minimizes $J$ in $\mathcal A$. Every minimizer has this form, and
\begin{equation}\label{hc:eq:exactmin}
 \min_{\mathcal A}J=3^{-\frac{3}{4}}D_0\int_0^U b(s)^3ds>0.
\end{equation}
Each $u_c$ is smooth and strictly increasing in the interior and gives a symphonic Hopf construction there.
\end{theorem}
\begin{proof}
For $X,Y\ge0$, $X^4+Y^4\ge4\,3^{-\frac{3}{4}}XY^3$, with equality exactly when $X=3^{-\frac{1}{4}}Y$. Apply this with $X=A^{\frac{1}{4}}|u'|$ and $Y=C^{\frac{1}{4}}b(u)$. For every $u\in\mathcal A$, the chain rule on interior intervals followed by the endpoint limits gives
\[
 J(u)\ge3^{-\frac{3}{4}}D_0\int_0^T|u'|b(u)^3dt
       \ge3^{-\frac{3}{4}}D_0\int_0^U b(s)^3ds.
\]
The simple zeros of $b$ make $\Psi:(0,U)\to\R$ a smooth increasing diffeomorphism. The endpoint hypotheses on $Z$ thus make \eqref{hc:eq:explicit} connect the two poles. Its derivative satisfies
\begin{equation}\label{hc:eq:firstorder}
 u'=\left(\frac C{3A}\right)^{\frac{1}{4}}b(u)
      =h\left(\frac q3\right)^{\frac{1}{4}}b(u)>0.
\end{equation}
Both inequalities are equalities, and changing variables to $u$ shows the energy is finite and equals \eqref{hc:eq:exactmin}. This proves existence and the exact infimum.

Conversely, equality for a minimizer forces the nonnegative Young defect to vanish almost everywhere and forces the primitive of $b^3$ composed with $u$ to be nondecreasing. Thus \eqref{hc:eq:firstorder} holds almost everywhere. A solution cannot depart from a pole at an interior point: on a compact interior interval its right side is locally Lipschitz in $u$, with the pole as an equilibrium. It therefore stays strictly inside, and integration gives \eqref{hc:eq:explicit} with one real constant. Finally, one may verify the second-order equation directly: \eqref{hc:eq:firstorder} gives $P=3^{-\frac{3}{4}}D_0b(u)^3$, hence $P'=Cb(u)^3b'(u)$. Corollary~\ref{hc:prop:full} completes the proof.
\end{proof}

\begin{remark}
This is a positive-infimum case with actual attainment, even though ordinary trace compactness is unavailable. The one-parameter family may concentrate at an endpoint as $c\to\pm\infty$. Thus failure of compactness of all minimizing sequences is compatible with the existence of individual minimizers. Neither positivity of the infimum alone nor a zero-capacity calculation settles attainment.
\end{remark}

\section{Smooth quaternionic constructions from one balance condition}\label{hc:sec:smooth}
The following statement handles all finite products at once. The construction with two source factors and the explicit additional scales below require only substitution into it.

\begin{theorem}[Balanced quaternionic slices]\label{short:thm:quaternion}
Fix $r\ge2$, $L=\frac{\pi}{2}$ and $R>0$. Let
$a_1=R\cos t$, $a_2=R\sin t$; for $i\ge3$, let $a_i>0$ be smooth on $[0,L]$ and even at both endpoints. Let $h>0$ be smooth and even there, with $h(0)=h(L)=R$. Put
\[
 g=h^2dt^2+\sum_{i=1}^r a_i^2g_{S^3},\qquad
 \Lambda=\sum_{i=1}^r a_i^{-2},\qquad
 \left(\prod_{i=1}^r a_i^2\right)\Lambda=K_0>0.
\]
Assume the last expression is constant. Then $g$ extends smoothly to
$M=S^7\times(S^3)^{r-2}$ (with the last product absent for $r=2$).
Let $S^4$ carry $G=du^2+b(u)^2g_{S^3}$, $0\le u\le U$, with $b>0$ inside and smooth pole expansions
\[
 b(u)=uB_-(u^2),\qquad b(U-v)=vB_+(v^2),\qquad B_-(0)=B_+(0)=1.
\]
For $f(x_1,\ldots,x_r)=x_1\cdots x_r$ in the unit quaternions, the minimizers of $J$ are precisely the functions
\[
 u_c(t)=\Psi^{-1}\!\left(\int_{\frac{\pi}{4}}^t h(\tau)\sqrt{\Lambda(\tau)}d\tau+c\right),
 \qquad c\in\mathbb R,\qquad
 J(u_c)=K_0^{\frac{3}{2}}\int_0^U b(s)^3ds,
\]
where $\Psi$ is defined in Theorem~\ref{hc:thm:calibration}.
Each gives a globally smooth symphonic map $M\to S^4$, and every angular factor is active.
\end{theorem}
\begin{proof}
Each partial map is left and right multiplication by unit quaternions, hence an isometry of $S^3$. The harmonic tension on each scaled product slice is zero, and
\[
 df\,df^*=\Lambda\operatorname{id},\quad
 \sigma_f=\Lambda df,\quad \mathcal T(f)=0,\quad
 Q_i=Q_{ij}=3,\quad q=3\Lambda^2.
\]
Here $\Lambda$ is constant on a slice. Therefore
\[
 D=\frac{\rho}{h}q^{\frac{3}{4}}=3^{\frac{3}{4}}K_0^{\frac{3}{2}},\qquad
 h\left(\frac{q}{3}\right)^{\frac{1}{4}}=h\sqrt\Lambda.
\]
At the left and right ends the latter function has the expansions
$\frac{1}{t}+tA_-(t^2)$ and $\frac{1}{L-t}+(L-t)A_+((L-t)^2)$, respectively. Theorem~\ref{hc:thm:calibration} gives the formula for $u_c$, characterizes all minimizers and gives their common energy.

For the geometry, in the left normal radius $v=R\sin t$ the radial coefficient is $\frac{h(t)^2}{R^2\cos^2t}$, a smooth function of $v^2$ equal to one at zero. The noncollapsing scales are also smooth functions of $v^2$. The reflected argument proves smoothness at the other stratum. The stated expansions of $b$ give a smooth target.

Near the first target pole, $\frac{1}{b(u)}=\frac{1}{u}+uD_-(u^2)$. Integrating the first-order equation for $u$ yields
\[
 \log u+G_-(u^2)=\log t+H_-(t^2)+c_-.
\]
The inverse function theorem for the odd function $u\exp G_-(u^2)$ gives $u=tF_-(t^2)$ with $F_-(0)>0$. Similarly $U-u=(L-t)F_+((L-t)^2)$.
In source normal coordinates $z=a_2(t)x_2$, the target normal expression is
$\frac{u}{a_2}x_1zx_3\cdots x_r$. At the other end, with $z=a_1(t)x_1$, it is
$\frac{U-u}{a_1}zx_2\cdots x_r$. The scalar coefficients are smooth in the squared normal radius and quaternion multiplication is multilinear. Thus the maps are globally smooth, and their interior symphonic equation extends by continuity. The partial isometries show that every factor acts nontrivially.
\end{proof}
\begin{corollary}[Quaternionic Hopf construction for general radial geometries]\label{hc:thm:quaternion}
Let $L=\frac{\pi}{2}$, $R>0$, and let $h\in C^\infty([0,L])$ be positive, even at both endpoints, with $h(0)=h(L)=R$. Equip the source $S^7$, parametrized by $(R\cos t\,x,R\sin t\,y)$, $x,y\in S^3$, with
\begin{equation}\label{hc:eq:quat-source}
 g=h(t)^2dt^2+R^2\cos^2t\,g_{S^3}+R^2\sin^2t\,g_{S^3}.
\end{equation}
Let the target $S^4$ carry any smooth rotationally symmetric metric
\begin{equation}\label{hc:eq:quat-target}
 G=du^2+b(u)^2g_{S^3},\qquad 0\le u\le U,
\end{equation}
with $b>0$ inside and smooth pole expansions
\[
 b(u)=u\,B_-(u^2),\qquad b(U-r)=r\,B_+(r^2),\qquad B_-(0)=B_+(0)=1.
\]
The functions $B_\pm$ are smooth near zero. No relation between $h$ and $b$ is imposed. Identify $S^3$ with the unit quaternions and fix the angular map $f(x,y)=xy$. Then the solutions
\begin{equation}\label{hc:eq:quat-solution}
 u_c(t)=\Psi^{-1}\left(\int_{\frac{\pi}{4}}^{t}\frac{h(\tau)}{R\sin\tau\cos\tau}d\tau+c\right),\qquad c\in\R,
\end{equation}
give globally smooth symphonic Hopf constructions $H_{f,c}: (S^7,g)\to(S^4,G)$. They are precisely the minimizers of $J$ for this angular map, with
\begin{equation}\label{hc:eq:quat-energy}
 J(u_c)=R^3\int_0^U b(s)^3ds.
\end{equation}
\end{corollary}
\begin{proof}
Take $r=2$ in Theorem~\ref{short:thm:quaternion}. Then
$\Lambda=(R^2\sin^2t\cos^2t)^{-1}$ and
$(a_1^2a_2^2)\Lambda=R^2$, so $K_0^{\frac{3}{2}}=R^3$. All geometric hypotheses coincide with those of that theorem, which gives the formula for $u_c$, its energy and the smoothness of the resulting map.
\end{proof}
For $R=1$, $h=1$, $U=\pi$, $b(u)=\sin u$, one has $\Psi(u)=\log\tan\left(\frac{u}{2}\right)$ and
\[
 u_c(t)=2\arctan(e^c\tan t),\qquad J(u_c)=\frac{4}{3}.
\]
At $c=0$ this is $u=2t$, a standard quadratic Hopf construction. The theorem also allows arbitrary smooth one- or multi-peaked target warps satisfying the pole conditions and arbitrary independent source functions $h$ satisfying the smooth-collapse conditions. Neither metric need be an induced ellipsoid metric. For example, the independent choices
\[
 h(t)=R\bigl(1+\delta\sin^2(2t)\bigr),\qquad
 b(u)=\sin u\bigl(1+\varepsilon\sin^2u\bigr),\qquad U=\pi,
\]
are allowed for any $\delta,\varepsilon>-1$. If $-1<\varepsilon<-\frac{1}{3}$, the target warping has two maxima separated by an interior local minimum, so this example lies beyond the one-peak hypothesis of Theorem~\ref{hc:thm:exist}.

\begin{corollary}[Smooth constructions with every angular factor active]\label{hc:cor:quaternion-multiple}
Let $r\ge3$, $L=\frac{\pi}{2}$, $R>0$, and $K>\frac{R^4}{4}$. For $4\le j\le r$ choose constants $\ell_j>0$, and put
\[
 S=\sum_{j=4}^r\ell_j^{-2},\qquad B=\prod_{j=4}^r\ell_j^2,
 \qquad p(t)=R^4\sin^2t\cos^2t.
\]
For $r=3$ use $S=0$ and $B=1$. Define
\begin{align}
 a_1&=R\cos t,\qquad a_2=R\sin t,\notag\\
 a_3^2&=\frac{K-p(t)}{R^2+S p(t)},\qquad a_j=\ell_j\quad(4\le j\le r).
 \label{hc:eq:multiple-warps}
\end{align}
Let $h>0$ be smooth and even at both endpoints, with $h(0)=h(L)=R$. The source metric
\[
 g=h(t)^2dt^2+\sum_{i=1}^ra_i(t)^2g_{S^3}
\]
extends smoothly to $M=S^7\times(S^3)^{r-2}$, collapsing the second factor at zero and the first factor at $L$. Let the target $(S^4,G)$ have any smooth rotationally symmetric metric satisfying the pole conditions of Corollary~\ref{hc:thm:quaternion}. Fix
\[
 f:(S^3)^r\to S^3,\qquad f(x_1,\ldots,x_r)=x_1x_2\cdots x_r,
 \qquad \Lambda(t)=\sum_{i=1}^ra_i(t)^{-2}.
\]
Then, for every $c\in\R$, the solution
\begin{equation}\label{hc:eq:multiple-solution}
 u_c(t)=\Psi^{-1}\left(\int_{\frac{\pi}{4}}^t h(\tau)\sqrt{\Lambda(\tau)}\,d\tau+c\right)
\end{equation}
gives a globally smooth symphonic Hopf construction $M\to S^4$. Every angular factor acts nontrivially. These functions exhaust the minimizers of $J$ for this fixed $f$, and
\begin{equation}\label{hc:eq:multiple-energy}
 J(u_c)=(BK)^{\frac{3}{2}}\int_0^U b(s)^3ds.
\end{equation}
The additional source warps obey \eqref{hc:eq:multiple-warps}; the assertion is not for arbitrary independent choices of all $a_i$.
\end{corollary}
\begin{proof}
Put $x=a_1^2$, $y=a_2^2$, $v=a_3^2$. Since $x+y=R^2$ and $xy=p$,
\[
 \left(\prod_i a_i^2\right)\Lambda
 =B[xy+v(x+y+Sxy)]=BK.
\]
The inequality $p\le \frac{R^4}{4}<K$ makes $a_3$ positive; its square is smooth and even at both endpoints. The remaining scales are positive constants. Thus all hypotheses of Theorem~\ref{short:thm:quaternion} hold with $K_0=BK$, giving every assertion.
\end{proof}
\begin{remark}[Finite products versus multiple joining parameters]\label{hc:rem:multiple}
Here the radial parameter is one-dimensional for every $r$. An $r$-fold topological join has $r-1$ radial parameters. In Corollary~\ref{hc:cor:quaternion-multiple}, every factor contributes to the angular map while the radial component depends only on $t$.
\end{remark}

\section{Zero infimum for finite products with linear collapse}\label{focus:sec:threshold}
This section concerns the reduced energy of the Hopf construction with one angular map on an arbitrary finite product. The relevant endpoint dimension is the sum of the dimensions of the factors collapsing there.

Let $r\ge2$, let $(P_i,g_i)$ be closed connected Riemannian manifolds of dimensions $m_i\ge1$, and let $f:\prod_iP_i\to(N,\gamma)$ be smooth, with $(N,\gamma)$ any closed Riemannian manifold. On the regular source put
\[
 g=h(t)^2dt^2+\sum_{i=1}^r\alpha_i(t)^2g_i,\qquad
 \rho=h\prod_i\alpha_i^{m_i},\qquad 0<t<T.
\]
All coefficients are smooth and positive inside, and $h$ has positive finite limits at both ends. For $e\in\{0,T\}$ write $\ell_0=t$, $\ell_T=T-t$. Choose nonempty sets $I_e\subset\{1,\ldots,r\}$ such that
\[
 \frac{\alpha_i}{\ell_e}\longrightarrow\sigma_{ei}>0\quad(i\in I_e),\qquad
 \alpha_j\longrightarrow a_{ej}>0\quad(j\notin I_e).
\]
Define $d_e=\sum_{i\in I_e}m_i$. Then $\rho\asymp\ell_e^{d_e}$ at that end.

The target regular metric is $k(s)^2ds^2+B(s)^2\gamma$, $-S<s<S$. Assume $k$ is positive and continuous on $[-S,S]$, $B$ is continuous there, positive inside and zero at both poles; both are smooth inside. Let $q_i,q_{ij}$ be the constant diagonal and mixed densities of $f$, or their normalized product-volume averages, as in \eqref{hc:eq:qt}. Set
\begin{align}
 Q_0(t)&=\sum_i\frac{q_i}{\alpha_i^4}
       +2\sum_{i<j}\frac{q_{ij}}{\alpha_i^2\alpha_j^2},\notag\\
 \mathcal H(\varphi)&=\frac14\int_0^T\rho\left[
       \frac{k(\varphi)^4}{h^4}|\varphi'|^4+Q_0(t)B(\varphi)^4\right]dt.
 \label{new:eq:hopf-weight}
\end{align}
This is $\frac{E(H_f)}{2\prod_i\Vol(P_i)}$. Let $\mathcal A$ consist of continuous functions valued in $[-S,S]$, locally $W^{1,4}$, connecting $-S$ to $S$ and having finite $\mathcal H$.

\begin{lemma}[Pairwise density constraints and the collapsing block]\label{new:density-cone}
For constant densities or normalized averages,
\begin{equation}\label{new:eq:density-cone}
 0\le q_{ij}\le\sqrt{q_iq_j}\qquad(i<j).
\end{equation}
Put
\begin{equation}\label{focus:eq:kappa}
 \kappa_e=\sum_{i\in I_e}\frac{q_i}{\sigma_{ei}^4}
       +2\sum_{\substack{i<j\\i,j\in I_e}}
           \frac{q_{ij}}{\sigma_{ei}^2\sigma_{ej}^2}.
\end{equation}
If $\kappa_e>0$, then $Q_0\sim\kappa_e\ell_e^{-4}$. If $\kappa_e=0$, then $Q_0$ is bounded near $e$; all densities involving a collapsing factor vanish.
\end{lemma}
\begin{proof}
In orthonormal frames write the partial differentials as $A_i$. The densities are $|A_iA_i^T|^2$ and $\operatorname{tr}(A_iA_i^TA_jA_j^T)$. Hilbert--Schmidt Cauchy--Schwarz gives the pointwise inequality; integration and Cauchy--Schwarz give the averaged version. The leading term of $Q_0$ is precisely \eqref{focus:eq:kappa}. If it vanishes, nonnegativity forces $q_i=0$ for every $i\in I_e$, and \eqref{new:eq:density-cone} forces $q_{ij}=0$ whenever either index collapses. Only bounded noncollapsing terms remain.
\end{proof}

\begin{theorem}[Complete zero-infimum criterion for finite products]\label{new:hopf-threshold}
For these data,
\begin{equation}\label{new:eq:full-threshold}
 \inf_{\mathcal A}\mathcal H=0
 \quad\Longleftrightarrow\quad
 \exists e\in\{0,T\}:\quad d_e>3
       \ \text{or}\ (d_e=3\text{ and }\kappa_e=0).
\end{equation}
A zero infimum is unattained and remains zero among admissible functions smooth in the interior and constant at the respective poles near both endpoints. Otherwise the infimum is finite and strictly positive. No angular symphonic equation is required for this energy assertion; positivity alone does not establish attainment.
\end{theorem}
\begin{proof}
A smooth transition on a compact interior interval has finite energy, so the admissible class is nonempty. It suffices to examine the left endpoint; the right endpoint follows by reflection.

\emph{An admissible transition and four separate estimates.}
Write $d=d_0$ and $I=I_0$. Fix $\eta\in C^\infty(\mathbb R,[0,1])$ nondecreasing, equal to zero on $(-\infty,\frac{1}{2}]$ and to one on $[1,\infty)$. For $0<\epsilon<\epsilon_0<\frac{T}{2}$, define
\begin{equation}\label{ep:eq:competitor}
 \varphi_\epsilon(t)=-S+2S\eta\left(\frac{t}{\epsilon}\right).
\end{equation}
It has the prescribed traces, is constant at the respective poles outside $[\frac{\epsilon}{2},\epsilon]$, and satisfies $|\varphi_\epsilon'|\le \frac{2S\|\eta'\|_\infty}{\epsilon}$. Both energy integrands vanish off that interval. For every fixed $\epsilon$ all coefficients are bounded on the transition interval, so $\varphi_\epsilon\in\mathcal A$.

Split $Q_0=Q_{00}+Q_{01}+Q_{11}$, where
\begin{align*}
 Q_{00}&=\sum_{i\in I}\frac{q_i}{\alpha_i^4}
       +2\sum_{\substack{i<j\\i,j\in I}}\frac{q_{ij}}{\alpha_i^2\alpha_j^2},\\
 Q_{01}&=2\sum_{\substack{i\in I\\j\notin I}}\frac{q_{ij}}{\alpha_i^2\alpha_j^2},\qquad
 Q_{11}=\sum_{i\notin I}\frac{q_i}{\alpha_i^4}
       +2\sum_{\substack{i<j\\i,j\notin I}}\frac{q_{ij}}{\alpha_i^2\alpha_j^2}.
\end{align*}
Here $q_{ji}=q_{ij}$, and empty sums are zero. Let $R_{00},R_{01},R_{11}$ be the corresponding sums with all denominators omitted. The collapse assumptions imply, after decreasing $\epsilon_0$,
\begin{equation}\label{ep:eq:weight-bounds}
 \rho\asymp t^d,\qquad Q_{00}\le CR_{00}t^{-4},\qquad
 Q_{01}\le CR_{01}t^{-2},\qquad Q_{11}\le CR_{11}.
\end{equation}
The positive functions $h,k$ are bounded above and away from zero on their closed intervals, and $B$ is bounded. Throughout the following estimates, $C$ may change but depends only on the fixed geometry, dimensions and cutoff, not on $\epsilon$. The four contributions to $4\mathcal H(\varphi_\epsilon)$ satisfy
\begin{align}
 E_{\rm rad}&:=\int_{\frac{\epsilon}{2}}^{\epsilon}
       \rho\frac{k(\varphi_\epsilon)^4}{h^4}|\varphi_\epsilon'|^4dt
       \le C\epsilon^{-4}\int_{\frac{\epsilon}{2}}^{\epsilon}t^d dt
       \le C\epsilon^{d-3},\label{ep:eq:radial}\\
 E_{00}&:=\int_{\frac{\epsilon}{2}}^{\epsilon}\rho Q_{00}B(\varphi_\epsilon)^4dt
       \le CR_{00}\int_{\frac{\epsilon}{2}}^{\epsilon}t^{d-4}dt
       \le CR_{00}\epsilon^{d-3},\label{ep:eq:collapsing}\\
 E_{01}&:=\int_{\frac{\epsilon}{2}}^{\epsilon}\rho Q_{01}B(\varphi_\epsilon)^4dt
       \le CR_{01}\int_{\frac{\epsilon}{2}}^{\epsilon}t^{d-2}dt
       \le CR_{01}\epsilon^{d-1},\label{ep:eq:mixed}\\
 E_{11}&:=\int_{\frac{\epsilon}{2}}^{\epsilon}\rho Q_{11}B(\varphi_\epsilon)^4dt
       \le CR_{11}\int_{\frac{\epsilon}{2}}^{\epsilon}t^d dt
       \le CR_{11}\epsilon^{d+1}.\label{ep:eq:noncollapsing}
\end{align}
Indeed $t=\epsilon v$ gives
$\int_{\frac{\epsilon}{2}}^{\epsilon}t^pdt=\epsilon^{p+1}\int_{\frac{1}{2}}^1v^pdv$ for every real $p$; the latter integral is finite, including $p=-1$. Thus no integrability at $t=0$ is assumed in these four estimates. Combining them gives
\begin{equation}\label{ep:eq:combined}
 0\le\mathcal H(\varphi_\epsilon)\le
 C\bigl[(1+R_{00})\epsilon^{d-3}
          +R_{01}\epsilon^{d-1}+R_{11}\epsilon^{d+1}\bigr].
\end{equation}
For $d>3$ every exponent is positive, which proves zero infimum. At the right endpoint use $S-2S\eta\left(\frac{T-t}{\epsilon}\right)$, the coordinate $T-t$, and the set $I_T$; the same estimate holds with $d_T$ and its corresponding density sums. This is one argument for arbitrary finite products, including simultaneous collapse of several factors.

\emph{The critical dimension with vanishing collapsing density.}
When $d_0=3$, the radial bound above is only $O(1)$ and does not prove concentration. If also $\kappa_0=0$, Lemma~\ref{new:density-cone} makes $Q_0$ bounded near zero. Choose $\chi\in C^\infty(\mathbb R,[0,1])$ nondecreasing, zero on $(-\infty,0]$ and one on $[1,\infty)$. For $0<\delta<\epsilon<\epsilon_0$, set
\begin{equation}\label{ep:eq:log-cutoff-function}
 \Lambda=\log\left(\frac{\epsilon}{\delta}\right),\qquad
 \varphi_{\delta,\epsilon}(t)=-S+2S\chi\!\left(\frac{\log\left(\frac{t}{\delta}\right)}{\Lambda}\right)
 \quad(t>0).
\end{equation}
Extend its value at zero by $-S$. This competitor is constant on $[0,\delta]$ and $[\epsilon,T]$, and $|\varphi_{\delta,\epsilon}'|\le \frac{C}{t\Lambda}$. Since $\rho\asymp t^3$, the two surviving energy terms obey
\begin{align}
 \int_\delta^\epsilon\rho\frac{k(\varphi_{\delta,\epsilon})^4}{h^4}
              |\varphi_{\delta,\epsilon}'|^4dt
 &\le C\Lambda^{-4}\int_\delta^\epsilon t^{-1}dt
   =C\Lambda^{-3},\label{ep:eq:log-radial}\\
 \int_\delta^\epsilon\rho Q_0B(\varphi_{\delta,\epsilon})^4dt
 &\le C\int_\delta^\epsilon t^3dt\le C\epsilon^4.
 \label{ep:eq:log-angular}
\end{align}
Consequently
\begin{equation}\label{new:eq:log-cutoff}
 \mathcal H(\varphi_{\delta,\epsilon})\le
 C\{\left[\log\!\left(\frac{\epsilon}{\delta}\right)\right]^{-3}+\epsilon^4\}.
\end{equation}
Taking $\epsilon=e^{-n}$ and $\delta=e^{-n^2}$ for sufficiently large integers $n$ proves zero infimum. Reflection gives the right-endpoint construction. The cutoffs are smooth and constant at their poles near both endpoints. No angular symphonic equation or pointwise constancy of the original densities was used; normalized angular averages suffice for this energy argument.

Conversely, suppose the condition in \eqref{new:eq:full-threshold} fails. Then $d_e\le3$, and $\kappa_e>0$ at every critical end $d_e=3$. If $\mathcal H(\varphi_j)\to0$, target arclength $u_j=K(\varphi_j)$, $K'=k$, has derivative energy $\int\left(\frac{\rho}{h^4}\right)|u_j'|^4\to0$. Boundedness, weighted H\"older on interior compact intervals and a subsequence at one interior point give local uniform convergence to a constant, so $\varphi_j\to s_*$ locally uniformly.

At an end with $d_e<3$, the integral of $\left(\frac{\rho}{h^4}\right)^{-\frac{1}{3}}$ converges. H\"older from that endpoint forces $s_*$ to equal its prescribed pole. If both ends are subcritical, this is already a contradiction. If a critical end exists, its $\kappa_e>0$ implies $Q_0(t)>0$ throughout the interior. Vanishing potential energy forces $B(s_*)=0$, so $s_*$ is one pole. Take the endpoint prescribing the opposite pole. The subcritical case is impossible by the same trace argument. In the critical case the lemma gives, on a fixed endpoint neighborhood,
\[
 \frac{\rho}{h}Q_0^{\frac{3}{4}}\ge c_0>0.
\]
Young's inequality gives the local lower bound
\[
 \mathcal H(\varphi_j)\ge3^{-\frac{3}{4}}c_0
       \int k(\varphi_j)|\varphi_j'|B(\varphi_j)^3dt.
\]
Apply the chain rule to a primitive of $kB^3$ and pass to the continuous endpoint trace. The limit inferior of the right-hand side is bounded below by
$3^{-\frac{3}{4}}c_0\int_{-S}^S k(s)B(s)^3ds$, contradicting vanishing energy.

Finally, zero energy forces $\varphi'=0$ almost everywhere on every interior compact interval, contrary to the distinct endpoint traces. This proves nonattainment when the infimum is zero.
\end{proof}

\begin{corollary}[Two-factor collapse]\label{focus:cor:two-threshold}
For $r=2$, $I_0=\{2\}$ and $I_T=\{1\}$, the criterion becomes
\[
 \inf\mathcal H=0\quad\Longleftrightarrow\quad
 \exists i\in\{1,2\}:\ m_i>3\ \text{or}\ (m_i=3\text{ and }q_i=0).
\]
\end{corollary}
\begin{proof}
Here $d_0=m_2$, $d_T=m_1$, and the corresponding $\kappa_e$ is a positive multiple of $q_2$ or $q_1$. Apply Theorem~\ref{new:hopf-threshold}.
\end{proof}
If several factors collapse at once, their dimensions must be added and their entire collapsing block enters $\kappa_e$. The scalar criterion remains valid even when the completed source has singular links. A claim about smooth global maps still requires smooth geometric completions; on such completions, the pole-constant competitors are globally smooth.

\section{Hopf constructions from arbitrary products to ellipsoids}\label{sec:hopf}
\subsection{One angular map on an arbitrary finite product}
Fix $r\ge2$ and arbitrary connected Riemannian manifolds $(P_i,g_i)$ of positive dimensions $m_i$. In this subsection compactness is not required for the local equation. Put
\begin{equation}\label{focus:eq:product}
 \mathcal P=\prod_{i=1}^r P_i,\qquad
 g=h(t)^2dt^2+g_t,\qquad g_t=\sum_{i=1}^r\alpha_i(t)^2g_i,
 \quad 0<t<T,
\end{equation}
where all coefficients are smooth and positive. Define $\rho=h\prod_i\alpha_i^{m_i}$.
The target, independently of these source data, is the ellipsoid
\begin{equation}\label{eq:concise:047}
 Y=Q^n(c,d)=\{(w,z)\in\mathbb R^n\times\mathbb R:
                   \frac{|w|^2}{c^2}+\frac{z^2}{d^2}=1\},\qquad c,d>0.
\end{equation}
With $L=\frac{\pi}{2}$, its regular metric is
\begin{equation}\label{eq:hopf-target}
 G=k(s)^2ds^2+c^2\cos^2s\,g_{S^{n-1}},\qquad
 k(s)^2=c^2\sin^2s+d^2\cos^2s,\quad -L<s<L.
\end{equation}
For a single smooth map $f:\mathcal P\to S^{n-1}$, possibly depending on every factor, set
\begin{equation}\label{eq:hopf-map}
 H_f(t,x_1,\ldots,x_r)
 =\bigl(c\cos\varphi(t)f(x_1,\ldots,x_r),d\sin\varphi(t)\bigr).
\end{equation}
This is the Hopf construction. Unlike the separated join $\Phi=(c\cos\varphi f_1,d\sin\varphi f_2)$, it uses one map into one angular target. Its source need not be an ellipsoid, and the $P_i$ need not be spheres or have constant curvature.

Use the diagonal and mixed densities $Q_i,Q_{ij}$ in \eqref{hc:eq:qt}, now with the round angular target. Then
\begin{equation}\label{eq:Q}
 Q(t,z)=|f^*g_{S^{n-1}}|_{g_t}^2
       =\sum_i\frac{Q_i(z)}{\alpha_i^4}
        +2\sum_{i<j}\frac{Q_{ij}(z)}{\alpha_i^2\alpha_j^2}.
\end{equation}
All mixed blocks are included. Writing $S_f=f^*g_{S^{n-1}}$, we obtain
\begin{equation}\label{eq:hopf-norm}
 |H_f^*G|_g^2=\frac{k(\varphi)^4}{h^4}(\varphi')^4
                        +c^4\cos^4\varphi\,Q(t,z).
\end{equation}
Define $L_t$ using the full slice metric, including its mixed raised pullback entries:
\begin{equation}\label{eq:Lt}
 L_tv=\frac1{\sqrt{\det g_t}}\partial_A
       (\sqrt{\det g_t}(S_f)^{AB}\partial_Bv),\qquad
 (S_f)^{AB}=(g_t)^{AC}(g_t)^{BD}(S_f)_{CD}.
\end{equation}

\begin{corollary}[Hopf criterion for arbitrary products]\label{thm:hopf}\label{lem:hopf}
For $\varphi\in C^2((0,T),(-L,L))$, the map \eqref{eq:hopf-map} is symphonic on its regular region if and only if
\begin{equation}\label{eq:hopf-angular}
 L_tf+Qf=0
\end{equation}
and
\begin{equation}\label{eq:hopf-radial}
 \frac1\rho\left(\rho\frac{k^4}{h^4}(\varphi')^3\right)'
 -\frac{k^3k'}{h^4}(\varphi')^4
 +c^4Q(t,z)\cos^3\varphi\sin\varphi=0
\end{equation}
hold pointwise for every $(t,z)$. The angular condition is exactly symphonicity of $f:(\mathcal P,g_t)\to S^{n-1}$ for every $t$.
\end{corollary}
\begin{proof}
Apply Corollary~\ref{new:master} with $H_s=c^2\cos^2s\,g_{S^{n-1}}$. Its angular pullback norm is $c^4\cos^4s\,Q$, whose negative quarter derivative is $c^4\cos^3s\sin s\,Q$. The sphere equation \eqref{eq:sphere} identifies the angular component. No source compactness, spherical geometry, or restriction to two source factors enters this substitution.
\end{proof}

On an interval with $0<|\varphi|<L$, Theorem~\ref{new:rank} determines exactly which combinations of the $\frac{r(r+1)}{2}$ densities must be constant. All individual densities are forced to be constant if the corresponding weight vectors span that full space on a detected interval. This span condition is automatic for the two ellipsoidal source warps below, but is an additional condition for general finite products. Symphonicity of the separate partial maps alone remains insufficient.

\subsection{Reduced energy and connecting solutions on finite products}
For the energy and global minimization statements in this subsection assume that the $P_i$ are closed. Denote constant densities or their normalized product-volume averages by $q_i,q_{ij}$, and put
\begin{align}
 V(t)&=c^4\left(\sum_i\frac{q_i}{\alpha_i^4}
                 +2\sum_{i<j}\frac{q_{ij}}{\alpha_i^2\alpha_j^2}\right),\notag\\
 J_H(\varphi)&=\frac14\int_0^T\rho\left[
       \frac{k(\varphi)^4}{h^4}(\varphi')^4+V(t)\cos^4\varphi\right]dt.
 \label{eq:JH}
\end{align}
Then $E(H_f)=2(\prod_i\Vol(P_i))J_H(\varphi)$, including when the energy is infinite. Averaging is legitimate for this identity, but the full equation still requires \eqref{eq:hopf-angular}--\eqref{eq:hopf-radial} pointwise.

\begin{corollary}[Connecting solutions for arbitrary closed angular factors]\label{focus:cor:ellipsoid-exist}
Suppose $f:(\mathcal P,g_t)\to S^{n-1}$ is symphonic and $Q(t,z)=q(t)>0$ is independent of $z$ for every $t$. If
\[
 \int_0^T\left(\frac{\rho}{h^4}\right)^{-\frac{1}{3}}dt<\infty,
\]
the polar connecting class $\varphi(0)=-L$, $\varphi(T)=L$ admits a minimizer that is smooth and strictly increasing on $(0,T)$. Its Hopf map is symphonic on the open cylinder. Under the endpoint and smooth-completion hypotheses of Proposition~\ref{hc:prop:global}, it extends to a globally weakly symphonic map. Alternatively, if the constant-$D$ and divergent-integral hypotheses of Theorem~\ref{hc:thm:calibration} hold, that theorem gives all minimizers explicitly without the finite-trace condition.
\end{corollary}
\begin{proof}
Set $u=\int_{-L}^{\varphi}k(s)ds$, $U=\int_{-L}^L k(s)ds$ and
$b(u)=c\cos\varphi(u)$. This warp has simple poles and one peak. The general product metric, density and normalized energy are exactly those in \eqref{hc:eq:source}--\eqref{hc:eq:energy}. Apply Theorem~\ref{hc:thm:exist}, Proposition~\ref{hc:prop:global}, or Theorem~\ref{hc:thm:calibration}, as specified.
\end{proof}

The pole limits send an entire endpoint stratum to $(0,\mp d)$ and respect any prescribed source collapse. For arbitrary closed factors the completion may have singular links. Accordingly, continuous pole limits and interior symphonicity do not by themselves assert a smooth global map between smooth closed manifolds.

\begin{proposition}[Local existence and uniqueness for the radial Hopf equation]\label{prop:local}
Let $\mathcal P$ and $g_t$ be as in \eqref{focus:eq:product}, allowing noncompact factors. Fix a smooth $f:\mathcal P\to S^{n-1}$ that is symphonic for every $g_t$, and suppose $Q(t,z)=q(t)$ is independent of $z$. Given $t_0\in(0,T)$, $s_0\in(-L,L)$ and $v_0\ne0$, the radial equation \eqref{eq:hopf-radial} has a smooth solution $\varphi:J\to(-L,L)$ on an open interval $J$ containing $t_0$, with
\[
 \varphi(t_0)=s_0,\qquad \varphi'(t_0)=v_0.
\]
The solution is unique as a germ at $t_0$. The associated map $H_f$ is symphonic on $J\times\mathcal P$.
\end{proposition}
\begin{proof}
Fix $z_0\in\mathcal P$. Evaluation at $z_0$ gives $q(t)=|f^*g_{S^{n-1}}|_{g_t}^2(z_0)$, a smooth nonnegative function. Let $D=\frac{\rho'}{\rho}-\frac{4h'}{h}$ and $V=c^4q$. Where $v\ne0$, \eqref{eq:hopf-radial} becomes
\[
 \varphi''=-\tfrac13D\varphi'-\frac{k'(\varphi)}{k(\varphi)}(\varphi')^2
 -\frac{h^4V\cos^3\varphi\sin\varphi}{3k(\varphi)^4(\varphi')^2}.
\]
Writing the displayed right-hand side as $G(t,s,v)$ gives the first-order system $(s',v')=(v,G(t,s,v))$. Its vector field is smooth on $0<t<T$, $|s|<L$, $v\ne0$, since $h,k$ are positive there. The Picard--Lindel\"of theorem gives a unique local solution through $(s_0,v_0)$ at $t_0$. Shrinking $J$ keeps its velocity nonzero, so the solved equation is equivalent to \eqref{eq:hopf-radial} on $J$. The hypothesis on $f$ supplies \eqref{eq:hopf-angular}, and Corollary~\ref{thm:hopf} gives the assertion about $H_f$.
\end{proof}

\begin{corollary}[An active nonspherical product example]\label{focus:cor:tori}
For each $i$, let $P_i$ be the quotient of $\mathbb R^{m_i}$ by the lattice $(2\pi\mathbb Z)^{m_i}$, equipped with any constant flat metric $g_i$, where $m_i\ge2$. Choose nonzero integer covectors $n_i$ and set
\[
 f(\theta_1,\ldots,\theta_r)=\exp\!\left(\mathrm i\sum_i n_i\cdot\theta_i\right)\in S^1.
\]
For arbitrary smooth positive $h,\alpha_i$ on $[0,T]$ and every target ellipsoid $Q^2(c,d)$, the radial equation admits a solution connecting the two poles. The associated symphonic Hopf map extends continuously from the open cylinder to the compact cylinder, sends its boundary components to opposite poles, and depends on every angular factor.
\end{corollary}
\begin{proof}
Write $\nu_i^2=|n_i|_{g_i^{-1}}^2>0$ and $\lambda(t)=\sum_i\frac{\nu_i^2}{\alpha_i^2}$. The circle phase has parallel differential on each flat product slice, so $\tau(f)=0$, $\sigma_f=\lambda df$ and $\mathcal T(f)=0$. Moreover $Q_i=\nu_i^4$, $Q_{ij}=\nu_i^2\nu_j^2$, and $Q=\lambda^2$ is constant in the angles. The source weight is positive on the closed interval, so the trace integral is finite. Corollary~\ref{focus:cor:ellipsoid-exist} applies. Each $n_i\ne0$ makes its factor active. This example concerns a cylinder with boundary; it makes no smooth closed-source collapse assertion for torus links.
\end{proof}

\subsection{Two sphere factors as a source specialization}\label{focus:subsec:spheres}
Only in the remainder of this section take $r=2$, $P_i=S^{m_i}$ with unit round metrics, $T=L$,
\[
 \alpha_1=a\cos t,\quad\alpha_2=b\sin t,\quad
 h^2=a^2\sin^2t+b^2\cos^2t.
\]
The completed source is then $Q^{m_1+m_2+1}(a,b)$. In its ambient coordinates,
\begin{equation}\label{eq:concise:010}
 H_f(a\cos t\,x,b\sin t\,y)
       =(c\cos\varphi(t)f(x,y),d\sin\varphi(t)).
\end{equation}
For consistency with the earlier ellipsoidal join notation, in the remainder of this specialization remove the constant source scale from the volume weight and energy:
\[
 \rho=h\cos^{m_1}t\sin^{m_2}t,\qquad
 J_H=\frac14\int_0^L\rho
 \left[\frac{k(\varphi)^4}{h^4}(\varphi')^4+V\cos^4\varphi\right]dt.
\]
This constant renormalization leaves all equations and minimizers unchanged. Now
$E(H_f)=2C_0J_H$, with $C_0$ as in \eqref{eq:concise:005}, and
\begin{equation}\label{eq:mixed-densities}
 Q=\frac{Q_1}{a^4\cos^4t}+\frac{Q_2}{b^4\sin^4t}
       +\frac{2Q_{12}}{a^2b^2\cos^2t\sin^2t}.
\end{equation}
We give a counterexample to slice-wise compatibility, followed by a rigidity criterion and a smooth example.
\begin{proposition}\label{prop:slices}
There is a map $f:S^1\times S^1\to S^3$ with both families of slices symphonic and all three densities constant, while its full angular equation fails for every $t\in I$.
\end{proposition}
\begin{proof}
Take $r,s>0$ with $r^2+s^2=1$, $r\ne s$, and put
\begin{equation}\label{eq:concise:048}
f(e^{i\theta},e^{i\eta})=(re^{i(\theta+\eta)},se^{i(\theta-\eta)}),\qquad
\delta=r^2-s^2\ne0.
\end{equation}
Since $|f_\theta|=|f_\eta|=1$ and $f_{\theta\theta}=f_{\eta\eta}=-f$, both slice families are unit speed great circles, hence symphonic. On the other hand,
\begin{equation}\label{eq:concise:049}
S_f=d\theta^2+2\delta\,d\theta d\eta+d\eta^2,
\qquad Q_1=Q_2=1,\quad Q_{12}=\delta^2.
\end{equation}
Put $\mu_1=(a^2\cos^2t)^{-1}$ and $\mu_2=(b^2\sin^2t)^{-1}$. Then
\begin{equation}\label{eq:concise:050}
L_tf=\mu_1^2f_{\theta\theta}+2\delta\mu_1\mu_2f_{\theta\eta}
+\mu_2^2f_{\eta\eta},\qquad
Q=\mu_1^2+\mu_2^2+2\delta^2\mu_1\mu_2.
\end{equation}
Consequently
\begin{equation}\label{eq:concise:051}
L_tf+Qf=2\delta\mu_1\mu_2(f_{\theta\eta}+\delta f),
\qquad f_{\theta\eta}+\delta f=(-2rs^2e^{i(\theta+\eta)},2r^2se^{i(\theta-\eta)}),
\end{equation}
is nowhere zero. Slice symphonicity therefore cannot replace the full angular equation.
\end{proof}
\begin{corollary}[Rigidity at an interior nonequatorial latitude]\label{prop:hopf-rigidity}
For the two-sphere source in Section~\ref{focus:subsec:spheres}, assume $H_f$ is symphonic and there exists $t_0\in I$ with $0<|\varphi(t_0)|<L$. Then $Q_1,Q_2,Q_{12}$ are all independent of the angles.
At a constant equatorial or polar latitude the coefficient of $Q$ in the radial equation vanishes.
\end{corollary}
\begin{proof}
Continuity gives an open interval about $t_0$ on which $\cos^3\varphi\sin\varphi\ne0$. The three weights span $\mathbb R^3$, by the polynomial independence check following Theorem~\ref{new:rank}. That theorem forces all three densities to be constant.
\end{proof}
\subsection{Endpoint concentration}
\begin{corollary}\label{prop:concentration}
For arbitrary smooth $f:\mathcal P\to S^{n-1}$, define $J_H$ using constant or averaged densities. Let
\begin{equation}\label{eq:concise:053}
\begin{split}
\mathcal A_{\mathrm{pol}}=\{\varphi\in C([0,L])\cap W^{1,4}_{\loc}(I):
&-L\le\varphi\le L,\\
&\varphi(0)=-L,\ \varphi(L)=L,\ J_H(\varphi)<\infty\}.
\end{split}
\end{equation}
If $m_2>3$ or $m_1>3$, then $\inf_{\mathcal A_{\mathrm{pol}}}J_H=0$ and the infimum is unattained. It remains zero even among smooth competitors constant near both endpoints; all corresponding competitors are globally smooth maps.
\end{corollary}
\begin{proof}
In Theorem~\ref{new:hopf-threshold}, take $\alpha_1=a\cos t$, $\alpha_2=b\sin t$, $B(s)=c\cos s$, $T=S=L$, and the ellipsoidal $h,k$. Its energy is $a^{m_1}b^{m_2}J_H$. The asserted dimension condition gives zero infimum and nonattainment. More explicitly, \eqref{ep:eq:combined} with $I_0=\{2\}$ gives
\begin{equation}\label{ep:eq:ellipsoid-bound}
 J_H(\varphi_\epsilon)\le C\bigl[(1+q_2)\epsilon^{m_2-3}
                   +q_{12}\epsilon^{m_2-1}+q_1\epsilon^{m_2+1}\bigr].
\end{equation}
The four integrals in \eqref{ep:eq:radial}--\eqref{ep:eq:noncollapsing} become respectively the radial, second-factor, mixed and first-factor terms. At the other endpoint interchange the two factor indices. The pole-constant smooth competitors give globally smooth maps on the ellipsoid.
\end{proof}

\begin{remark}\label{rem:scope}
Corollary \ref{prop:concentration} rules out only absolute minimizers in $\mathcal A_{\mathrm{pol}}$, not nonzero energy critical points or other variational constructions. The comparison does not require the angular equation. At an end with $m_i=3$, the leading power in the cutoff estimate is zero. If $q_1,q_2>0$ and both $m_1,m_2\le3$, Corollary~\ref{prop:hopf-threshold} gives a positive lower bound. If the other collapsing dimension exceeds three, a minimizing sequence can concentrate at that end.
\end{remark}

\subsection{The exact dimension threshold for two positive factor densities}
\begin{corollary}[A necessary and sufficient condition for zero infimum]\label{prop:hopf-threshold}
Assume $q_1,q_2>0$ and $q_{12}\ge0$, where each $q$ is either a constant density or its angular average. Set
\begin{equation}\label{eq:hopf-V-threshold}
 V(t)=c^4\left[\frac{q_1}{a^4\cos^4t}+\frac{q_2}{b^4\sin^4t}
 +\frac{2q_{12}}{a^2b^2\cos^2t\sin^2t}\right].
\end{equation}
In the polar connecting class $\mathcal A_{\mathrm{pol}}$ in \eqref{eq:concise:053},
\begin{equation}\label{eq:sharp-threshold}
 \inf J_H=0\quad\Longleftrightarrow\quad\max(m_1,m_2)>3.
\end{equation}
The zero infimum is unattained. If $m_1,m_2\le3$, then
\begin{equation}\label{eq:hopf-positive-bound}
 J_H(\varphi)\ge\frac{c_*}{3^{\frac{3}{4}}}\int_{-L}^{L}k(s)\cos^3s\,ds>0,
 \qquad c_*:=\inf_{t\in I}\frac{\rho(t)}{h(t)}V(t)^{\frac{3}{4}}>0.
\end{equation}
This bound holds for every admissible function, including nonmonotone ones. It does not itself ensure attainment. Nor is equality asserted for every choice of angular data and semiaxes. Indeed, put $D(t)=\frac{\rho}{h}V^{\frac{3}{4}}$. For a smooth admissible function $\varphi$ with $\varphi'>0$ and $-L<\varphi<L$ on $I$, equality throughout the lower-bound argument requires $D(t)=c_*$ on $I$ and
\[
 \frac{k(\varphi)^4}{h^4}\varphi'^4=\frac13V(t)\cos^4\varphi.
\]
These are additional conditions. The calibrated families established earlier verify their own equality cases.
\end{corollary}
\begin{proof}
Theorem~\ref{new:hopf-threshold} under the same substitution gives the equivalence and nonattainment because $q_1,q_2>0$. For the explicit bound, $c_*>0$ follows from
\begin{equation}\label{eq:threshold-asymptotics}
 \frac{\rho}{h}V^{\frac{3}{4}}\asymp t^{m_2-3}\quad(t\downarrow0),\qquad
 \frac{\rho}{h}V^{\frac{3}{4}}\asymp(L-t)^{m_1-3}\quad(t\uparrow L).
\end{equation}
The Young inequality used in the proof of that theorem, followed by the chain rule for $G'(s)=k(s)\cos^3s$, gives \eqref{eq:hopf-positive-bound}. No monotonicity is needed.
\end{proof}
\noindent Positivity of both $q_1,q_2$ is essential to this endpoint comparison. If either vanishes, \eqref{eq:threshold-asymptotics} may change; Theorem~\ref{new:hopf-threshold} gives the corrected criterion in that case.

\begin{proposition}[The classical quadratic Hopf map]\label{prop:quadratic-hopf}
For $m_1=m_2=1$, $n=2$, and $a=b=c=d=1$, take
\begin{equation}\label{eq:concise:056}
f(e^{i\theta},e^{i\eta})=e^{i(\theta+\eta)},\qquad \varphi(t)=2t-L.
\end{equation}
Then $H_f:S^3\to S^2$ is globally smooth and symphonic.
\end{proposition}
\begin{proof}
Here $Q_1=Q_2=Q_{12}=1$. Put $\mu_1=\cos^{-2}t$ and $\mu_2=\sin^{-2}t$. Then
\begin{equation}\label{eq:concise:057}
L_tf=-(\mu_1+\mu_2)^2f,\qquad
Q=(\mu_1+\mu_2)^2=(\sin t\cos t)^{-4},
\end{equation}
so the full angular equation holds. Moreover $h=k=1$, $\rho=\sin t\cos t$, and $\varphi'=2$. The two radial terms are
\begin{equation}\label{eq:concise:058}
8(\cot t-\tan t),\qquad
Q\cos^3\varphi\sin\varphi=\frac{8(\sin^2t-\cos^2t)}{\sin t\cos t},
\end{equation}
and cancel. With $z_1=\cos t\,e^{i\theta}$ and $z_2=\sin t\,e^{i\eta}$, the map is
\begin{equation}
H_f(z_1,z_2)=\bigl(2z_1z_2,|z_2|^2-|z_1|^2\bigr).
\label{eq:quadratic-hopf}
\end{equation}
This polynomial map has squared norm $(|z_1|^2+|z_2|^2)^2=1$ and is smooth on all of $S^3$. Its symphonic divergence vanishes everywhere by continuity from the dense regular region.
\end{proof}

\section*{AI assistance}
The work \cite{Cao} was done withou any AI tools, Three years later, I used ChatGpt to improve it.   The author used generative AI tools, including OpenAI Codex, in developing and revising mathematical arguments, editing the manuscript, checking selected formulas and references, and preparing Lean proofs for selected statements. The author takes responsibility for the mathematical claims, references and final text.

\end{document}